\documentclass[10pt,a4paper,oneside,reqno]{amsart}

\usepackage{setspace}
\usepackage[totalwidth=13.5cm,totalheight=22cm]{geometry}
\usepackage[T1]{fontenc}
\usepackage[dvips]{graphicx}
\usepackage{epsfig}
\usepackage{amsmath,amssymb,amscd,amsthm}
\usepackage{subfigure}
\usepackage[latin1]{inputenc}
\usepackage{latexsym}
\usepackage{graphics}
\usepackage{colortbl} 
\usepackage{color}
\usepackage{ae}
\usepackage{enumitem}
\usepackage[all]{xy}
\usepackage{scalerel}
\usepackage{tikz}
\usepackage{cancel}
\usepackage{mathrsfs}  
\usepackage[colorlinks=true,pagebackref=true]{hyperref}
\usepackage[normalem]{ulem}
\usepackage{nicefrac}
\usepackage{xfrac}

\makeatletter
\newtheorem*{rep@theorem}{\rep@title}
\newcommand{\newreptheorem}[2]{%
\newenvironment{rep#1}[1]{%
 \def\rep@title{#2 \ref{##1}}%
 \begin{rep@theorem}}%
 {\end{rep@theorem}}}
\makeatother

\makeatletter
\newtheorem*{rep@cor}{\rep@title}
\newcommand{\newrepcor}[2]{%
\newenvironment{rep#1}[1]{%
 \def\rep@title{#2 \ref{##1}}%
 \begin{rep@cor}}%
 {\end{rep@cor}}}
\makeatother

\makeatletter
\newtheorem*{rep@prop}{\rep@title}
\newcommand{\newrepprop}[2]{%
\newenvironment{rep#1}[1]{%
 \def\rep@title{#2 \ref{##1}}%
 \begin{rep@prop}}%
 {\end{rep@prop}}}
\makeatother

\newtheorem{cor}{Corollary}[section]

\newrepcor{cor}{Corollary}
\newreptheorem{theorem}{Theorem}
\newrepprop{prop}{Proposition}

\theoremstyle{definition}
\newtheorem{defi}[cor]{Definition}
\theoremstyle{remark}
\newtheorem{remark}[cor]{Remark}
\newtheorem*{remark*}{Remark}

\newtheorem*{notation*}{Notation}

\newlist{steps}{enumerate}{1}
\setlist[steps, 1]{itemsep=8pt,leftmargin=0cm,itemindent=.5cm,labelwidth=\itemindent,labelsep=0cm,align=left,label = \textbf{\emph{Step \arabic*}:\,}}

\newcommand{\R}{{\mathbb R}}
\newcommand{\Z}{{\mathbb Z}}

\newcommand{\Hyp}{\mathbb{H}}
\newcommand{\Sph}{\mathbb{S}}
\newcommand{\AdS}{\mathbb{A}\mathrm{d}\mathbb{S}}
\newcommand{\E}{\mathbb{E}}
\newcommand{\M}{\mathbb{M}}
\newcommand{\HP}{\mathbb{HP}}

\newcommand{\RP}{\mathbb{R}\mathrm{P}}
\newcommand{\SP}{\mathbb{R}\mathrm{S}}

\newcommand{\GL}{\mathrm{GL}}

\newcommand{\Isom}{\mathrm{Isom}}

\renewcommand{\O}{\mathrm{O}}

\newcommand{\p}[1]{\ensuremath{\boldsymbol{#1^+} }}
\newcommand{\m}[1]{\ensuremath{\boldsymbol{#1^-} } }
\renewcommand{\l}[1]{\ensuremath{\boldsymbol{#1}} }

\begin{document}

\setcounter{secnumdepth}{2}
\setcounter{tocdepth}{1}

\title[Examples of geometric transition in low dimensions]{Examples of geometric transition in low dimensions}

\author[Andrea Seppi]{Andrea Seppi}
\address{Andrea Seppi: CNRS and Universit\'e Grenoble Alpes \newline 100 Rue des Math\'ematiques, 38610 Gi\`eres, France. } \email{andrea.seppi@univ-grenoble-alpes.fr}

\thanks{The author is member of the national research group GNSAGA}

\begin{abstract}
The purpose of this note is to discuss examples of geometric transition from hyperbolic structures to half-pipe and Anti-de Sitter structures in dimensions two, three and four. As a warm-up, explicit examples of transition to Euclidean and spherical structures are presented. No new results appear here; nor an exhaustive treatment is aimed. On the other hand, details of some elementary computations are provided to explain certain techniques involved. This note, and in particular the last section, can also serve as an introduction to the ideas behind the four-dimensional construction of \cite{rioloseppi}.
\end{abstract}

\maketitle

\vspace{-0.8cm}

\section*{Introduction}

The idea of \emph{geometric transition}, or \emph{transitional geometry}, goes back to Klein (\cite{MR1509812}, see also \cite{AP}) and to the observation that Euclidean geometry can be seen as a \emph{transition} between spherical and hyperbolic geometry. In a modern language, a geometric transition consists in a deformation of geometric structures (by which we mean $(G,X)$-structures, see for instance Thurtson's famous notes \cite{thurstonnotes}) which \emph{degenerate}, but admit a limit consisting in a different type of geometric structure on the same manifold. For instance, there are examples of hyperbolic structures which degenerate to a single point, and admit a Euclidean structure as a limit. Viceversa, one talks of \emph{regeneration} when a Euclidean structure can be deformed into hyperbolic or spherical structures, see 
\cite{P98} where the existence of such phenomenon in dimension three is proved under certain cohomological conditions.

Recently, in his PhD thesis \cite{danciger} Danciger introduced a different phenomenon, which is a transition between hyperbolic geometry and its Lorentzian counterpart, namely Anti-de Sitter geometry (that is, Lorentzian geometry of constant negative sectional curvature, whose study on the mathematical point of view has been pioneered by Mess \cite{Mess} and has largely developed in recent times). Here the transitional geometry between the two is the so-called \emph{half-pipe} geometry, which was in fact known since before under the name of co-Minkowski geometry, being naturally a dual geometry for Minkowski space. Half-pipe geometry is the ``right'' geometry to describe the degeneration of hyperbolic structures towards a codimension one object, a pheonomenon which has been largely studied and played an important role in the Orbifold Theorem (see among others \cite{Hthesis,CHK,BLP,MR2140265,P13}). In \cite{dancigertransition} (see also \cite{dancigerideal}) the existence of phenomena of regeneration from half-pipe structures in dimension three was proven, under the hypothesis of a cohomological condition in the spirit of \cite{P98}. The recent paper \cite{rioloseppi} provided the first examples of geometric transition in dimension four, on a certain class of cusped, finite-volume manifolds.

The purpose of this note is to describe rather explicitly several examples of geometric transition from hyperbolic to half-pipe and Anti-de Sitter geometry, first in dimension two and three, which is done in Section  \ref{sec:HP}, and then to outline the four-dimensional examples of \cite{rioloseppi}, which is the content of Section \ref{sec:4d}. But before that, in Section \ref{sec:eucl} certain examples of transition from hyperbolic to Euclidean and spherical geometry are provided as a sort of warm-up, again in dimension two and three. In fact, all  these examples in low dimension start from the same base point: in dimension two, from deforming a complete finite-volume hyperbolic structure on a punctured torus, so as to make it collapse in one case to a single point, in the other to a circle; in dimension three, we start instead from a  complete finite-volume hyperbolic structure on the Borromean rings complement, which is obtained from an ideal right-angled octahedron, and we deform it to obtain a collapse to a single point or to a codimension one totally geodesic surface.
Although mostly elementary, we provide some of the details of these examples in a rather precise fashion, with the purpose of highlighting the methods involved.
As already stated, in Section \ref{sec:4d} we give a brief outline of the four-dimensional construction of the recent work \cite{rioloseppi}, which is based on a deforming polytope introduced in \cite{KS} and \cite{MR}, and we rely on the previous lower dimensional examples to help the four-dimensional intuition.

It is also worth mentioning the work that has been done recently to study more general possibilities of geometric transition, see for instance \cite{CDW} and \cite{surveyseppifillastre} for constructions inside real projective geometry, and \cite{trettel_thesis} for new ideas which leave the ambient projective setting. We believe it would be interesting to provide new classes of examples for all these phenomena.

\subsection*{Acknowledgements}

I would like to thank Stefano Riolo and Pierre Will for their interest and for several conversations.

\section{From hyperbolic structures...}

Let us start by recalling some notions of hyperbolic geometry. This section will then provide two fundamental examples of hyperbolic manifolds.

\subsection{Models of hyperbolic space} The \emph{hyperbolic space} of dimension $n$ is defined, in the \emph{hyperboloid model}, as:
$$\Hyp^n=\{x=(x_0,\ldots,x_n)\in\R^{n+1}\,:\,q_{1,n}(x)=-1,\,x_0>0\}~,$$
where we denote by $q_{1,n}$ the quadratic form
$$q_{1,n}(x)=-x_0^2+x_1^2+\ldots+x_n^2$$
and more in general
\begin{equation} \label{eq:quadratic forms}
q_{a,n+1-a}(x)=-x_0^2\ldots -x_{a-1}^2+x_a^2+\ldots+x_n^2~.
\end{equation}
This definition highlights the fact that $\Hyp^n$ is essentially the analogue in Minkowski space of the $n$-sphere $\Sph^n$, which is just
$$\Sph^n=\{x\in\R^{n+1}\,:\,q_{0,n+1}(x)=1\}~.$$
Endowed with the Riemannian metric which is the restriction of the standard bilinear form of signature $(1,n)$ on $\R^{n+1}$, whose associated quadratic form is $q_{1,n}$, the hyperbolic space is the unique complete, simply connected Riemannian manifold of constant sectional curvature -1 up to isometries. Its group of isometries $\Isom(\Hyp^n)$ is isomorphic to $\O_+(1,n)$, namely, the index two subgroup of $\O(1,n)$ which preserves $\Hyp^n\subset\R^{n+1}$.

Totally geodesic subspaces of $\Hyp^n$ of dimension $k$ are obtained as the intersections of $\Hyp^n\subset\R^{n+1}$ with $(k+1)$--dimensional linear subspaces of $\R^{n+1}$ (when such intersection is non-empty). In particular geodesics are the intersections with linear planes, and hyperplanes of $\Hyp^n$ are the intersection of $\Hyp^n$ with linear hyperplanes of $\R^{n+1}$. To every such hyperplane of the form $P\cap \Hyp^n$, one can associate a reflection, which is the unique element of $\Isom(\Hyp^n)$ which fixes pointwise $P$ and acts as minus the identity on its orthogonal complement with respect of the standard bilinear form of signature $(1,n)$ of $\R^{n+1}$.

There is also a well-known \emph{projective model}, or \emph{Klein model}, of $\Hyp^n$, which is the image of the obvious inclusion $\Hyp^n\hookrightarrow\RP^n$ sending $x$ to $[x]$. That is, we consider $\Hyp^n$ as the domain
\begin{equation} \label{eq:domain in proj space}
\{[x]\in\RP^n\,:\,q_{1,n}(x)<0\}\subset\RP^n
\end{equation}
and in this model the isometry group of $\Hyp^n$ corresponds to the group of projective transformations which preserve such domain. By the above discussion, totally geodesic subspaces in this model are the intersections of the domain \eqref{eq:domain in proj space} with projective subspaces.  The projective model is particularly helpful to visualise the \emph{boundary at infinity} $\partial \Hyp^n$ of $\Hyp^n$, namely the boundary of $\Hyp^n$ inside $\RP^n$:
$$\partial\Hyp^n=\{[x]\in\RP^n\,:\,q_{1,n}(x)=0\}~,$$
which is homeomorphic to a sphere of dimension $n-1$.

The hyperbolic space $\Hyp^n$ is the local model for \emph{hyperbolic manifolds}, meaning that every Riemannian manifold $\mathcal M$ of constant sectional curvature -1 can be equivalently seen as the datum of a $(\Isom(\Hyp^n),\Hyp^n)$--structure on $\mathcal M$: namely, an atlas for $\mathcal M$ with values in $\Hyp^n$ and such that the transition functions, where defined, are the restrictions of isometries of $\Hyp^n$. The purpose of this note is to exhibit several explicit examples of deforming hyperbolic structures, which will collapse \emph{as hyperbolic structures} but admit a limit which is in fact another kind of geometric structure on $\mathcal M$. 

The common strategy to build such geometric structures is by constructing a \emph{convex hyperbolic polytope}, namely a region in $\Hyp^n$ obtained as the intersection of a finite family of half-spaces (that is, regions bounded by hyperplanes), and identifying the codimension one faces in pairs by means of certain isometries of $\Hyp^n$. We shall see here two fundamental examples, namely an ideal quadrilateral (in dimension 2) and an ideal octahedron (in dimension 3). Starting from these two examples, in Sections \ref{sec:eucl} and \ref{sec:HP} we will then explore phenomena of geometric transition, towards Euclidean/spherical geometry and half-pipe/Anti-de Sitter geometry respectively.

\subsection{Ideal quadrilateral and the punctured torus} \label{sec:quadrilateral}

\begin{figure}[htb]
\includegraphics[height=4cm]{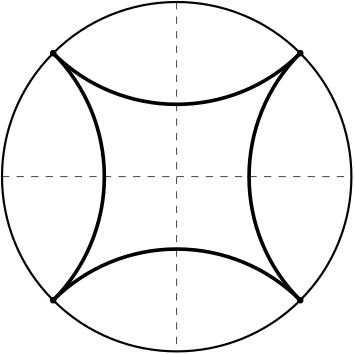}
\caption{\footnotesize A regular ideal quadrilateral in $\Hyp^2$ and its axes of symmetry, in the Poincar\'e (conformal) model.}\label{fig:quadrilateral}
\end{figure}

Let us start by a simple example in dimension 2. The polygon we consider here is a \emph{regular ideal quadrilateral}, namely a quadrilateral $\mathcal Q$ bounded by four geodesics which do not intersect in $\Hyp^2$ but whose closures in $\partial\Hyp^2$ intersect in four different points at infinity. Moreover, $\mathcal Q$ will be maximally symmetric, meaning that it will be preserved by a dihedral group of order 8 generated by a rotation of angle $\pi/2$ and a reflection. See Figure \ref{fig:quadrilateral}. 
The standard identification of the sides which produces (topologically) a torus out of a quadrilateral can be realised by isometries, thus obtaining a complete, finite-area hyperbolic metric on a punctured torus. Despite being very elementary, we provide here some details of this construction since the following examples will be computed exactly along the same lines.

\begin{remark} \label{rmk linear forms}
In this and the following examples, we will express half-planes in $\Hyp^2$ (and in general  in $\Hyp^n$) as elements of the dual projective sphere $\SP^{n,*}=\R^{n+1,*}/\R_{+}$. Namely, we consider the space of non-trivial linear forms $\alpha:\R^{n+1}\to\R$, and associate to every such linear form the half-space $\{x\in\R^{n+1}\,:\,\alpha(x)\leq 0\}$. Two such linear forms give the same half-space if and only if they are positive multiples of one another. The condition that the hyperplane defined as the kernel of $\alpha$ intersects $\Hyp^n$, and therefore $\alpha$ defines a honest half-space of $\Hyp^n$ (i.e. the intersection is not the empty set or the full $\Hyp^n$) is easily expressed by the fact that $q_{1,n}^*(\alpha)<0$, where $q_{1,n}^*$ is the quadratic form  of signature $(1,n)$ on $\R^{n+1,*}$, in standard form with respect to the dual basis of the standard basis of $\R^{n+1}$.
\end{remark}

Back to the concrete example of a regular ideal quadrilateral, we can assume for simplicity of computations that the center of symmetry is the point $e_0=(1,0,0)$. Hence the configuration of planes in $\R^3$ will be invariant with respect to the Euclidean rotation of angle $\pi/2$ around the axis spanned by $e_0$. (In fact such Euclidean rotation also preserves the quadratic form $q_{1,n}$, and thus induces an element of $\Isom(\Hyp^3)$.)

\begin{notation*}
Hereafter we shall use the notation $[x_0:\ldots:x_n]$ to denote projective classes of elements of $\R^{n+1}$, expressed in the standard basis, whereas the notation $(y_0:\ldots:y_n)$ will denote projective classes in the dual vector space $\R^{n+1,*}$, in the coordinate system given by the dual standard basis.
\end{notation*}

We can moreover assume that the (ideal) vertices of the quadrilateral are the points of $\partial\Hyp^2$ defined by 
\begin{equation} \label{eq:vertices ideal quad}
[1:\pm \sqrt 2/2:\pm\sqrt 2/2]~.
\end{equation}
The four half-planes bounding $\mathcal Q$ thus have the form 
$(-1:\pm \sqrt 2:0)$ and $(-1:0:\pm \sqrt 2)$, as it can be easily checked that each linear form $\alpha$ in this projective classes vanishes on precisely two of the ideal vertices, and moreover $\alpha(e_0)<0$. The four sides of the quadrilateral are then identified by means of hyperbolic isometries: the first maps the ``left'' side to the ``right'' side, and is obtained by composing the reflection $(x_0,x_1,x_2)\mapsto (x_0,-x_1,x_2)$ with a reflection in one of these two sides; the second maps the ``top'' side to the ``bottom'' side and is again the composition of $(x_0,x_1,x_2)\mapsto (x_0,x_1,-x_2)$ with the reflection in one of the sides. Glueing the sides according to these instructions gives rise to a complete, finite-area hyperbolic structure on a punctured torus.

In Sections \ref{subsec:2d eucl} and \ref{subsec:3d eucl} we will see how this hyperbolic structure on the punctured torus can be deformed to non-complete structures having cone singularities at the puncture. In one case these deformed structures will collapse to a single point, and will be continued to Euclidean and spherical structures after suitable rescaling: in the second case, the deformed structures collapse to a closed geodesic and, suitably rescaled, produce a geometric transition to half-pipe and Anti-de Sitter geometry.

\subsection{The three-dimensional ideal octahedron} \label{sec: octahedron}

\begin{figure}[htb]
\includegraphics[height=5cm]{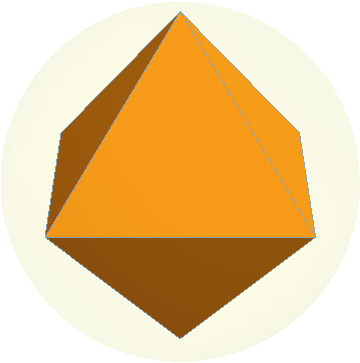}
\caption[The cuboctahedron.]{\footnotesize The ideal octahedron in $\Hyp^3$, seen in the projective model.}\label{fig:octahedron}
\end{figure}

Let us now rise to dimension three. We consider the \emph{right-angled ideal octahedron} $\mathcal O$, which is a maximally symmetric octahedron in $\Hyp^3$ with totally geodesic faces and vertices at infinity, and such that all dihedral angles are $\pi/2$. See Figure \ref{fig:octahedron}. Similarly to the case of the ideal quadrilateral, we can assume that the center of symmetry is $e_0$ and place the vertices at the points 
\begin{equation} \label{eq: vertices oct}
[1:\pm \sqrt 2/2:\pm\sqrt 2/2:0]\qquad \text{and}\qquad[1:0:0:\pm 1]~.
\end{equation} The first four vertices are actually the vertices of a regular ideal quadrilateral in the totally geodesic plane defined by $x_3=0$. Symmetries along planes are again easily written as the reflections in $\R^4$ which change sign to $x_1$, $x_2$ or $x_3$.

When the combinatorics of a polytope becomes more intricate (as we shall see later), it becomes simpler to describe the polytope as the intersection of half-spaces instead of listing vertices. For the octahedron the defining (projective classes of) linear forms, seen as elements of the dual projective sphere $\SP^3$ (see Remark \ref{rmk linear forms}) are collected in Table \ref{table:oct ideal}.

\begin{table}[htb]
\begin{eqnarray*}
\left( -1:-\sqrt{2}:0:-1 \right) , & &
\left( -1:-\sqrt{2}:0:+1 \right),\\
\left( -1:0:-\sqrt{2}:+1 \right),& &
\left( -1:0:-\sqrt{2}:-1 \right),\\
\left( -1:+\sqrt{2}:0:-1 \right),& &
\left( -1:+\sqrt{2}:0:+1 \right),\\
\left( -1:0:+\sqrt{2}:+1 \right),& &
\left( -1:0:+\sqrt{2}:-1 \right).\\
\end{eqnarray*}
\caption{\footnotesize Linear forms defining the right-angled ideal octahedron, expressed as elements of $\SP^{3,*}$. The sign of the last coordinates detects whether the corresponding face lies in the upper or lower half-space, as in Figure \ref{fig:octahedron}. One easily checks that $e_0$ is in the interior of all the corresponding half-spaces, hence in the interior of the ideal octahedron, and that faces meet in groups of four in the ideal vertices listed in \eqref{eq: vertices oct}. All dihedral angles are right.}\label{table:oct ideal}
\end{table}

It is known that the faces of the right-angled ideal octahedron can be paired by isometries so as to obtained a complete, finite-volume hyperbolic manifold. In fact different faces can be paired by means of an orientation-preserving isometry, using an even number of reflections either in a plane of symmetry or in a face. For instance in \cite[Section 3.3]{thurstonnotes} a construction which produces a hyperbolic manifold homeomorphic to the Whitehead link complement is described. There are actually many possible non-homeomorphic (hence non-isometric, by Mostow rigidity) complete hyperbolic manifolds which can be obtained, see \cite{MR2484431} for a complete discussion. In Sections \ref{subsec:2d HP} and \ref{subsec:3d HP} we will explore two deformations of the ideal octahedron which can be used to produce three-dimensional examples of geometric transition to Euclidean/spherical or to half-pipe/Anti-de Sitter structures.

\section{...to Euclidean and spherical structures}\label{sec:eucl}

Here we will describe two examples of geometric transition from hyperbolic to Euclidean and spherical structures, which deform the examples of Sections \ref{sec:quadrilateral} and \ref{sec: octahedron}.

\subsection{A ``classical'' geometric transition}
In order to describe and understand geometric transition from hyperbolic to Euclidean and spherical structures, we first need to explain how Euclidean geometry can be seen as a ``rescaled limit'' of both hyperbolic and spherical geometry. For this purpose, let us embed Euclidean space $\E^n$ into $\R^{n+1}$ by $(x_1,\ldots,x_n)\mapsto (1,x_1,\ldots,x_n)$, hence
with image the horizontal plane at height 1. Let us observe that such embedding is isometric if $\R^{n+1}$ is endowed either with the standard (Euclidean) positive definite quadratic form $q_{0,n+1}$ or with the (Minkowski) quadratic form $q_{1,n}$ which is used to define the hyperboloid model of $\Hyp^n$.

\begin{figure}[htb]
\includegraphics[height=6cm]{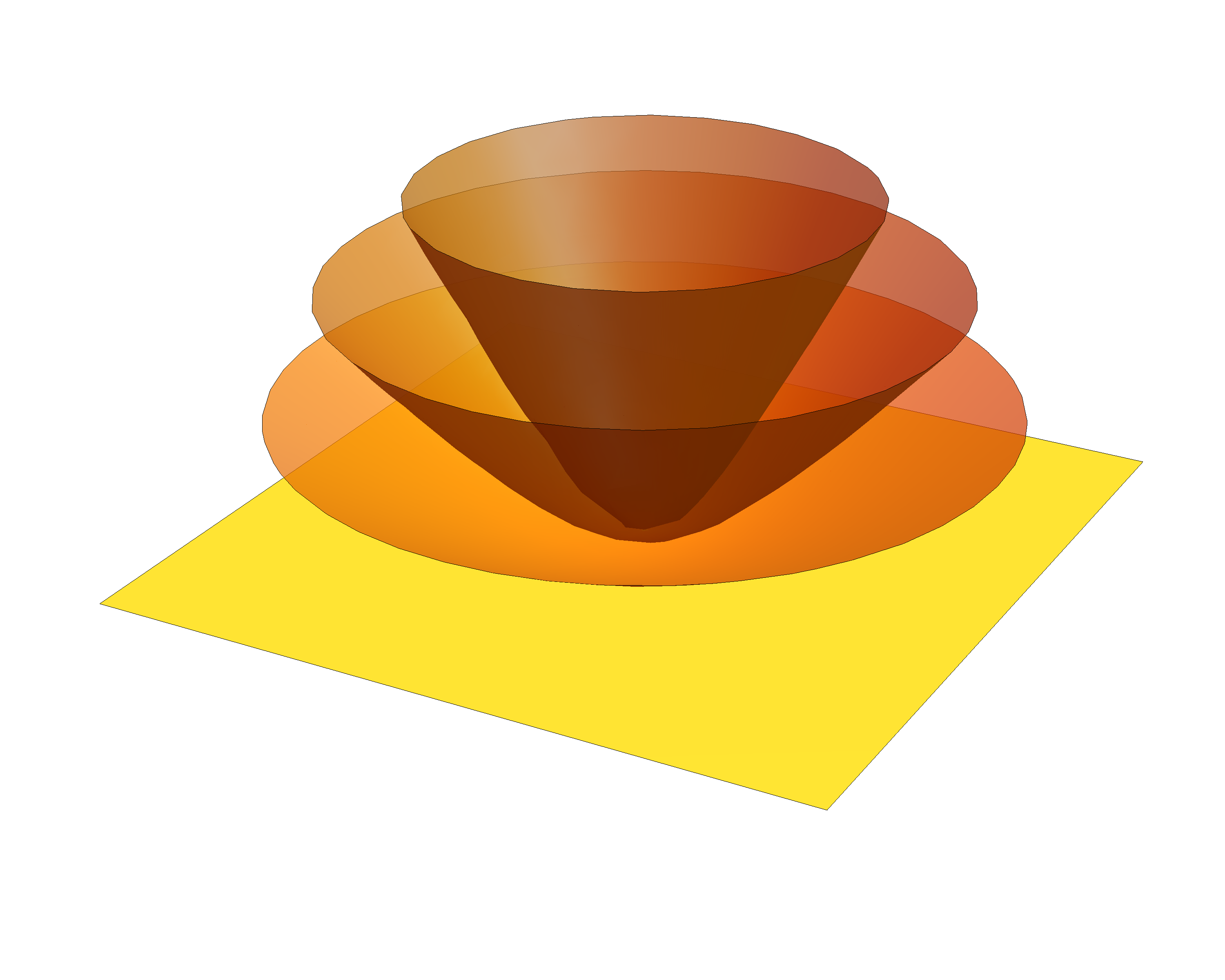}
\vspace{-1cm}
\caption{\footnotesize The Hausdorff limit of $\gamma_t(\Hyp^n)$ in $\R^{n+1}$ is the horizontal Euclidean plane at height 1. This picture represents the situation for $n=2$.}\label{fig:hausdorff1}
\end{figure}

Let us now consider the linear transformations $\gamma_t$ of $\R^{n+1}$ which fix the point $e_0$ and stretch all the other directions $e_1,\ldots,e_n$ by a factor $1/t$, namely in standard coordinates
\begin{equation} \label{eq gammat}
\gamma_t=\mathrm{diag}\left(1,\frac{1}{t},\ldots,\frac{1}{t}\right)~.
\end{equation}
We will use these transformations $\gamma_t$ to ``zoom in'' hyperbolic space $\Hyp^n$ (or the $n$-sphere $\Sph^n$) around the point $e_0$. More formally, the Hausdorff limit of the subsets $\gamma_t(\Hyp^n)\subset \R^{n+1}$ as $t\to 0^+$ is precisely $\E^n$, embedded as the horizontal plane at height 1. See Figure \ref{fig:hausdorff1}. The same holds for $\gamma_t(\Sph^n)$, except that the Hausdorff limit in this case is $\pm\E^n$, namely the two horizontal planes at height $\pm 1$. Moreover, one can show that the group of isometries of Euclidean space $\E^n$ is the Chabauty limit of the groups $\gamma_t\Isom(\Hyp^n)\gamma_t^{-1}$ and $\gamma_t\Isom(\Sph^n)\gamma_t^{-1}$, seen as closed subgroups of $\GL(n+1,\R)$, where we recall that $\Isom(\Hyp^n)\cong \O_+(1,n)$ and $\Isom(\Sph^n)\cong \O(n+1)$. See \cite{cooperdanwien} and \cite[Chapter 4]{surveyseppifillastre} for more details.

In fact, hyperbolic/Euclidean/spherical metrics on a manifold $\mathcal M$, which are equivalent to the data of a $(\Isom(X),X)$-structure on $\mathcal M$ for $X=\Hyp^n/\E^n/\Sph^n$ respectively, can all be seen as a special kind of \emph{real projective structures}. These are atlases for $\mathcal M$ with values in $\RP^n$ and such that the transition functions are projective transformations. To see this, it is more convenient to work with the \emph{real projective sphere} 
$\SP^n=\R^{n+1}/\R_{>0}$, which is the 2-to-1 cover of $\RP^n$, so that we can see $\Hyp^n$, $\E^n$ and $\Sph^n$ simultaneously as subsets of $\SP^n$: $\Hyp^n$ is defined by the conditions $q_{1,n}(x)<0$ and $x_0>0$, $\E^n$ simply by $x_0>0$ (all these conditions are preserved by multiplication by a positive number) and $\Sph^n$ is the whole $\SP^n$. Moreover, the isometry groups of the three spaces all embed into the group of projective transformations of $\SP^n$, $\GL(n+1,\R)/\R_{>0}$. Thus in conclusion, $(\Isom(X),X)$-structures are special forms of real projective structures, regardless as $X=\Hyp^n,\E^n$ or $\Sph^n$.

The above discussion permits to give the following definition of geometric transition, which will be realised concretely in several examples in the rest of Section \ref{sec:eucl}.

\begin{defi} \label{defi transition eucl}
A \emph{geometric transition} on a manifold $\mathcal M$ from hyperbolic to spherical geometry, through Euclidean geometry, is a continuous path of real projective structures $\mathscr P_t$ on $M$, defined for $t\in (-\epsilon,\epsilon)$, which is conjugate to
\begin{itemize}
\item Hyperbolic structures for $t>0$;
\item Euclidean structures for $t=0$;
\item Spherical structures for $t<0$.
\end{itemize}
\end{defi}

By a continuous path of real projective structures we mean, in short, a family of atlases which vary continuously on $t$ and such that the transition functions also vary continuously.

\subsection{On a (punctured) torus} \label{subsec:2d eucl}

The example of Section \ref{sec:quadrilateral} can be easily deformed to give rise to a geometric transition as in Definition \ref{defi transition eucl}. To visualise this, let us deform the regular ideal quadrilateral to a (non-ideal) quadrilateral $\mathcal Q(t)$ with vertices in $\Hyp^2$, which is again regular in the sense that it is invariant under rotation by a dihedral group of order 8, and we shall again assume that the center of symmetry is $e_0$. The vertices of such a deformation can be easily expressed as 
$$p_i(t)=\exp_{e_0}(tv_i)=\cosh(t)e_0+\sinh(t)v_i$$
where $\exp$ is the exponential map of $\Hyp^2$ and $v_i$ ($i=1,2,3,4$) are the unit tangent vectors $(0,\pm\sqrt 2/2,\pm\sqrt 2/2)$ in $T_{e_0}\Hyp^2=\{x_0=0\}$. Indeed, when $t\to +\infty$, $[p_i(t)]$ converges in $\RP^2$ to the vertices of the regular ideal quadrilateral (see \eqref{eq:vertices ideal quad}).

\begin{figure}[htb]
\includegraphics[width=\textwidth]{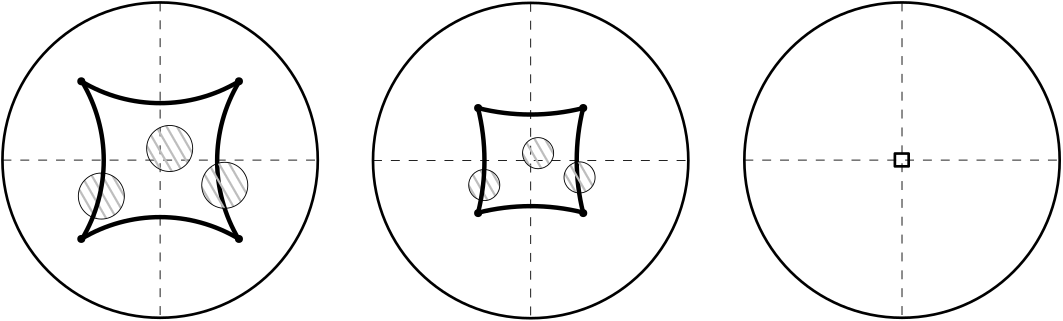}
\caption[The cuboctahedron.]{\footnotesize The collapsing family of quadrilaterals, together with a family of continuously varying charts for the associated hyperbolic structures.}\label{fig:quadrilateral2}
\end{figure}

By a simple computation, we then see that the \emph{rescaled limit} of $\mathcal Q(t)$ as $t\to 0^+$ is a Euclidean square. Namely, recalling that $\gamma_t$ is the linear map which stretches the directions $e_1$ and $e_2$ by a factor $1/t$, one has
$$\lim_{t\to 0^+}\gamma_t(p_t)=\lim_{t\to 0^+}\left(\cosh(t)e_0+\frac{\sinh(t)}{t}v_i\right)=e_0+v_i~,$$
that is, we obtain the four vertices $(1,\pm\sqrt 2/2,\pm\sqrt 2/2)$ of a square (of diagonal equal to $2$) in the Euclidean plane defined by $x_0=1$. It is also clear that the sides of the quadrilateral $\mathcal Q(t)$, which are infinite geodesics in $\Hyp^2$ obtained as intersections with linear planes in $\R^{3}$ and collapse to either $\{x_1=0\}$ or $\{x_2=0\}$ as $t=0$, in the rescaled limit converge to the sides of the square (which are again intersections of $\E^2$, seen as the horizontal plane $\{x_0=1\}$, with linear planes of $\R^3$). The angles at the vertices vary monotonically between $0$ and $\pi/2$.

There is a completely analogous deformation for spherical quadrilaterals. In fact, taking the exponential map of $\Sph^2$ and replacing $\sinh$ and $\cosh$ by $\sin$ and $\cos$, the very same argument shows that the same Euclidean square is the rescaled limit of a family of regular spherical quadrilaterals (which become non-convex at $t=\pi/2$, when the vertices reach the equator), whose angles vary monotonically between $\pi/2$ and $3\pi/2$ as $t$ varies in $(0,\pi)$. To help the intuition, we declare that this spherical deformation occurs for negative times (as in Definition \ref{defi transition eucl}), hence we call $\mathcal Q(-t)$ such a spherical quadrilateral.

\begin{remark}
There is of course nothing special of the quadrilateral above, nor of the dimension 2. One can apply the same construction to ``exponentiate'' any Euclidean polygon (or polytope, in higher dimension) to a family of hyperbolic/spherical polygons (or polytopes). But in higher dimension the geometric structures obtained by such construction will typically have a complicated singular locus, more than those we will provide, by a different method, in Section \ref{subsec:2d HP}.
\end{remark}

We are now ready to see that the complete, finite-area hyperbolic structure on the punctured torus given in Section \ref{sec:quadrilateral} can be deformed in the realm of real projective structures so as to obtain a geometric transition as in Definition \ref{defi transition eucl}. In fact, we can produce a path of hyperbolic structures $\mathscr H_t$ on the punctured torus, for $t\in(0,+\infty]$ again by identifying the sides of $\mathcal Q(t)$ in pairs: the identification of the left and right side is obtained as the composition of the reflection $(x_0,x_1,x_2)\mapsto (x_0,-x_1,x_2)$ in the vertical axis of symmetry with reflection in one of the sides, and similarly for the top and bottom sides. Such hyperbolic structures are non-complete, as the punctures can be reached in finite length. Their completion are obtained by just adding one point at the puncture, and is a Riemannian cone manifold with cone singularity approaching $2\pi$ as $t\to 0^+$ and $0$ as $t\to +\infty$. The same construction applies on the spherical side, thus giving rise to spherical structures $\mathscr S_t$ (here $-\pi<t<0$) on the punctured torus of cone angles larger than $2\pi$ at the puncture (as expected by the Gauss-Bonnet formula). 

Let $\mathscr P_t$ be the real projective structure obtained by:
\begin{itemize}
\item Post-composing the charts of the hyperbolic structures $\mathscr H_t$ with the projective transformation induced by $\gamma_t$, for $t>0$;
\item Post-composing the charts of the spherical structures $\mathscr S_t$ with the projective transformation induced by $\gamma_{|t|}$, for $t<0$;
\item Identifying the sides of the unit Euclidean square by horizontal and vertical translations, for $t=0$.
\end{itemize}
 The transition functions of $\mathscr P_t$ are the conjugate by $\gamma_t$ of those of $\mathscr H_t$ and $\mathscr S_t$. In fact one can easily produce charts for $\mathscr P_t$ which vary continuously in $t$, as in Figure \ref{fig:quadrilateral2} (which pictures the situation before rescaling). To verify that $\mathscr P_t$ is a geometric transition as in Definition \ref{defi transition eucl}, it only remains to show that the transition functions vary continuously. 

To see this, recall that transition functions (for all the cases: the hyperbolic, Euclidean and spherical structures) are obtained as the composition of two reflections, one along an axis of symmetry and another along one of the sides. Moreover, the order of composition of these do not depend on $t$, and is the same in all the three cases. Now, the reflections in the axis of symmetry have the form $(x_0,x_1,x_2)\mapsto(x_0,\mp x_1,\pm x_2)$ in all three geometries and for all $t$, and it is easily checked that such maps commute with $\gamma_t$. We claim that reflections in one of the sides of $\mathcal Q(t)$, when conjugate by $\gamma_t$, converge to the Euclidean reflection in the corresponding side of the Euclidean square. This will conclude the fact that the transition functions pass to the limit as expected when $t\to 0^\pm$, and therefore $\mathscr P_t$ define a geometric transition as in Definition \ref{defi transition eucl}. 

To check the claim, let $\alpha_t:\R^3\to\R$ be a linear form which vanishes on the linear plane containing a geodesic of $\mathcal Q(t)$ (hyperbolic for $t>0$, spherical for $t<0$) having unit norm for the dual quadratic form $q_{1,2}^*$ ($t>0$) or $q_{0,3}^*$ ($t<0$). Hence $\alpha_0(x)=\pm x_1$ or  $\alpha_0(x)=\pm x_2$, according to the side of the quadrilateral that had been chosen. The hyperbolic of spherical reflection matrix can then be written in the standard basis as 
\begin{equation}\label{eq:hyp reflection}
\rho_t=\mathrm{id}-2J_\pm \alpha_t\alpha_t^T~,
\end{equation}
where $J_\pm$ is the matrix $\mathrm{diag}(-1,1,1)$ if $t>0$ and the identity matrix if $t<0$, and $\alpha_t$ is now thought as a column vector whose entries are the coordinates of the unit linear form with respect to the standard dual basis $e_0^*,e_1^*,e_2^*$. Indeed, $\rho_t$ obviously preserves any vector in the kernel of $\alpha_t$, and one can check (using the assumption that $\alpha_t$ are unit linear forms) that it acts on its orthogonal complement for the hyperbolic/spherical metric, which is $J_\pm\alpha$, as minus the identity. A direct computation then shows that 
$$\lim_{t\to 0^\pm}\gamma_t\rho_t\gamma_t^{-1}=
\begin{pmatrix} 
\mathrm{id} &  0 \\
-2ca_0 & r 
\end{pmatrix}
 $$
 where $a_0\in\R^2$ is a vector such that $\alpha_0=(0,a_0)$, $c$ is a constant which is basically the derivative at $t=0$ of the first coordinate of $\alpha_t$, and $r=\mathrm{id}-2a_0^Ta_0$ is the matrix of the Euclidean reflection in $\R^2$ fixing $a_0^\perp$. In our concrete example, $c=\sqrt 2/2$ and $\lim_{t\to 0^\pm}\gamma_t\rho_t\gamma_t^{-1}$ acts on $\E^2$, seen as the horizontal plane at height 1, as $z\mapsto r(z)-2ca_0$, the reflection in the corresponding side of the Euclidean square of diagonal $2$.



\subsection{On the Borromean rings complement and the three-torus} \label{subsec:3d eucl}
We will now move to dimension three and exhibit a one-parameter family of polytopes $\mathcal O(t)$ which deform the octahedron (recall Figure \ref{fig:octahedron}) and collapse to a single point when $t\to 0^+$. Rescaling the deformation as in the previous example will permit to ``continue'' the deformation to a Euclidean and spherical polytope, and to construct a three-dimensional geometric transition. This example essentially appears in \cite[Chapter 3]{CHK}; some details of the proof will be omitted, since they follow the same lines as in the previous section.

\begin{figure}[htb]
\centering
\begin{minipage}[c]{.25\textwidth}
\centering
\includegraphics[scale=0.25]{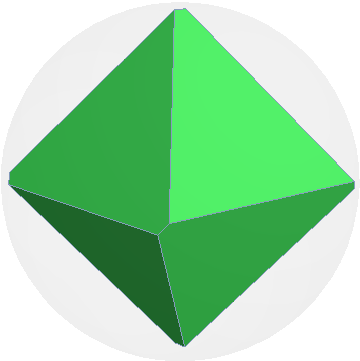}
\end{minipage}%
\begin{minipage}[c]{.25\textwidth}
\centering
\includegraphics[scale=0.25]{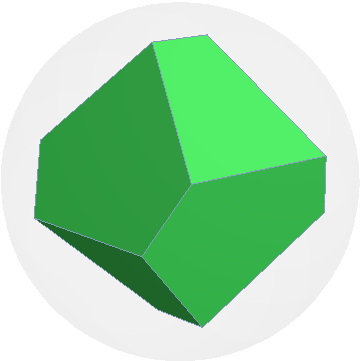}
\end{minipage}%
\begin{minipage}[c]{.25\textwidth}
\centering
\includegraphics[scale=0.25]{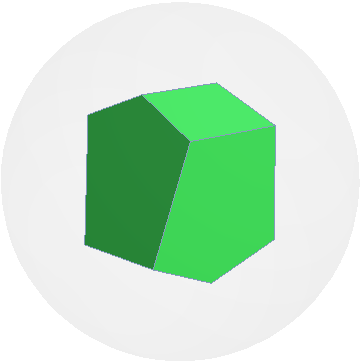}
\end{minipage}
\begin{minipage}[c]{.24\textwidth}
\centering
\includegraphics[scale=0.25]{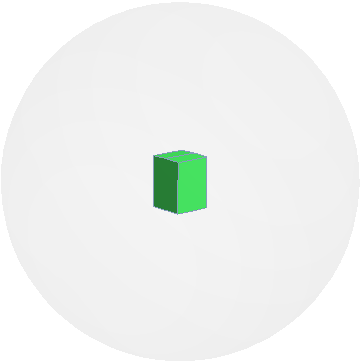}
\end{minipage}%
\caption{\footnotesize A movie of the deformation ${\mathcal O}(t)$ in an affine chart (Klein model of $\Hyp^3$), from the right-angled ideal octahedron to the collapse at a single point.}
\label{fig:collapse oct eucl}
\end{figure}

The idea of such a deformation is that the ideal vertices are moved towards the interior of $\Hyp^3$ and become edges of the new polytope with varying dihedral angle. These edges will thus give rise to cone singularities along circles in the geometric structures we obtain. All the other dihedral angles remain right-angled. When $t\to 0^+$, $\mathcal O(t)$ will degenerate to a single point. See Figure \ref{fig:collapse oct eucl}.


The list of half-spaces defining $\mathcal O(t)$ are collected in Table \ref{table:oct collapse1}, as elements of $\SP^{3,*}$ (Remark \ref{rmk linear forms}). For $t=1$ the list coincides with Table \ref{table:oct ideal}, namely the right-angled ideal octahedron. Let us briefly explain how to interpret this list. First, we clearly have that $e_0$ lies in the interior of each half-space (except for $t=0$), hence in the interior of $\mathcal O(t)$. The sign of the last coordinate determines whether the corresponding half-space is  ``upper'' (that is, defined by an inequality $x_3\geq f(x_1,x_2)$ in the Klein model of $\Hyp^3$ given by $x_0=1$) or ``lower'' ($x_3\leq f(x_1,x_2)$). 
There is a qualitative difference between the left column and the right column: the half-spaces in the left column converge to either the upper half-space $x_3\geq 0$ or the lower half-space $x_3\leq 0$ when $t\to 0$. Hence the corresponding faces tend to become horizontal, namely, contained in the totally geodesic plane $\{x_3=0\}$. On the other hand, the half-spaces in the right column converge to a half-space which is orthogonal to $\{x_3=0\}$ (and containing $e_0$ in its boundary). 

The parameter $t$ has the following geometric interpretation: the hyperbolic distance between the lower edge (which is the intersection of the faces given by the first and third element in the left column) and the upper edge (second and fourth on the left) equals $2\tanh(t)$.
One can check that the dihedral angle on these two edges is the same and approaches $\pi$ as $t\to 0^\pm$; the other four non-right dihedral angles approach $\pi/2$ instead. 

\begin{table}[htb]
\begin{eqnarray*}
\left( -t:-\sqrt{2}t^2:0:-1 \right) , & &
\left( -t:-\sqrt{2}:0:+t^2 \right),\\
\left( -t:0:-\sqrt{2}t^2:+1 \right) , & &
\left( -t:0:-\sqrt{2}:-t^2 \right),\\
\left( -t:+\sqrt{2}t^2:0:-1 \right) , & &
\left( -t:+\sqrt{2}:0:+t^2 \right),\\
\left( -t:0:+\sqrt{2}t^2:+1 \right) , & &
\left( -t:0:+\sqrt{2}:-t^2 \right).\\
\end{eqnarray*}
\caption{\footnotesize The half-spaces defining the deformation of the right-angled ideal octahedron, expressed as elements of $\SP^{3,*}$. When $t=1$ we recover the ideal octahedron, while for $t=0$ the polyhedron collapses to a single point.}\label{table:oct collapse1}
\end{table}

\begin{figure}[htb]
\includegraphics[height=3.5cm]{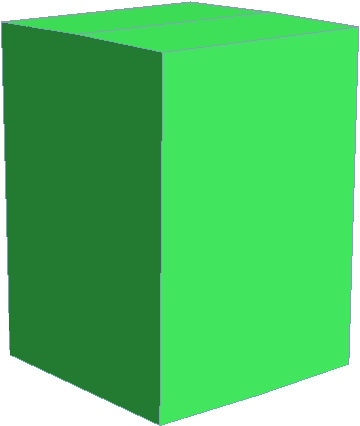}
\caption{\footnotesize The rescaled limit of the polytope $\mathcal O(t)$ as $t\to 0^\pm$ is a Euclidean parallelepiped.}\label{fig:parallelepiped rescaled}
\end{figure}

The limiting values of dihedral angles are not surprising, since as $t\to 0^\pm$ the geometry is approaching Euclidean geometry at an infinitesimal scale, and the polytope is approaching a Euclidean right-angled parallelepiped. This is seen more clearly when taking the limit of the rescaled polytope $\gamma_{t}\mathcal O(t)$. To compute this, it is sufficient to observe that applying $\gamma_{t}$ to $\mathcal O(t)$ just amounts to multiplying the first coordinate in Table \ref{table:oct collapse1} by $1/t$. Hence the limit as $t\to 0^+$ is the polytope in $\SP^3$ defined by the linear forms $(-1:0:0:-1)$ (which is the limit of each of the two ``lower' faces whose intersection is the lower edge), $(-1:0:0:1)$ (same for the upper edge),  $(-1:\pm \sqrt 2:0:0)$ and $(-1:0:\pm \sqrt 2:0)$ (the ``vertical'' faces). Hence we obtain the Euclidean parallelepiped in $\E^3$ (which we recall is obtained as $x_0=1$) defined by $-1\leq x_3\leq 1$, $-\sqrt 2/2\leq x_1\leq \sqrt 2/2$, $-\sqrt 2/2\leq x_2\leq \sqrt 2/2$. See Figure \ref{fig:parallelepiped rescaled}.

\begin{figure}[htb]
\centering
\begin{minipage}[c]{.25\textwidth}
\centering
\includegraphics[scale=0.2]{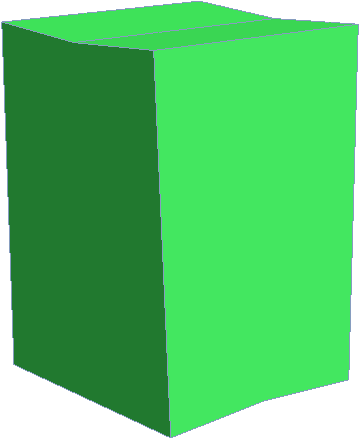}
\end{minipage}%
\begin{minipage}[c]{.25\textwidth}
\centering
\includegraphics[scale=0.2]{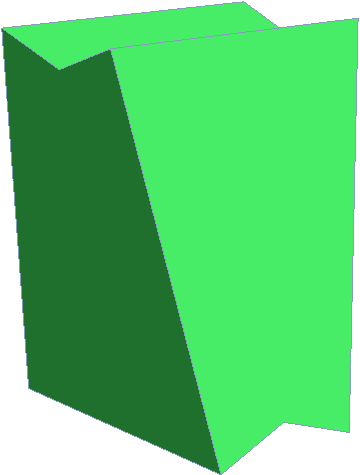}
\end{minipage}%
\begin{minipage}[c]{.25\textwidth}
\centering
\includegraphics[scale=0.2]{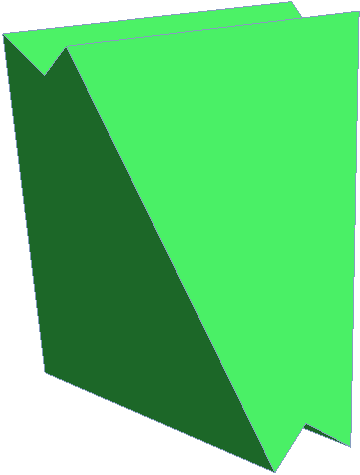}
\end{minipage}
\begin{minipage}[c]{.2\textwidth}
\centering
\includegraphics[scale=0.2]{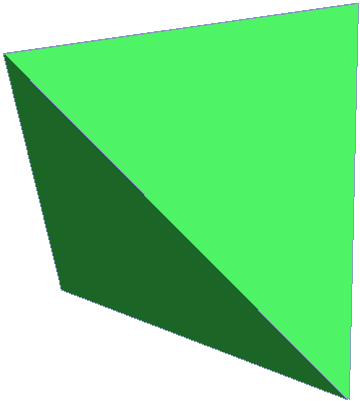}
\end{minipage}%
\caption{\footnotesize The deformation $\mathcal O(t)$ for negative times, rescaled by $\gamma_{|t|}$, which interpolates between the Euclidean parallelepiped of Figure \ref{fig:parallelepiped rescaled} and spherical polytopes with the same combinatorics. Some faces degenerate when $t=-1$.}
\label{fig:regen oct eucl}
\end{figure}

The deformation can be in fact continued for negative times to $\Sph^3$. A simple way to see this is to define $\mathcal O(-t)$ as the polytope in $\R^4\cap \Sph^3$ with the same combinatorics as $\mathcal O(t)$, whose faces lie in the hyperplanes as in the list of Table \ref{table:oct collapse1}, except that one must change the sign of the last three coordinates in the four linear forms of the right column. This choice assures that those planes bounding $\mathcal O(t)$ which are orthogonal in $\Hyp^3$ for $t>0$ remain orthogonal, now in $\Sph^3$, for $t<0$. However the spherical polytope obtained in this way is not convex, hence it cannot be defined simply as the intersection of half-spaces. See Figure \ref{fig:regen oct eucl}.

There are several ways to build a geometric transition from hyperbolic to Euclidean and spherical geometry out of the deformation $\mathcal O(t)$. Probably the simplest to visualise, at least on the hyperbolic side, is obtained by first identifying the faces in the right column of Table \ref{table:oct collapse1} in pairs (first and third, second and fourth), using as usual the symmetries which change sign to $x_1$ or $x_2$ and reflections in these planes, to obtain a manifold (homeomorphic to the product of a punctured torus and a closed interval) with totally geodesic boundary and corners along the upper and lower edges (which have now become circles). Then one doubles such manifold along the boundaries (that is, one glues two copies, with opposite orientations, using the identity on the boundary). This gives rise to the complement of three circles in the 3-torus, which is also homeomorphic to the complement of a link, the Borromean rings, in the 3-sphere. These hyperbolic structures have cone singularities along each component of the link, and the cone angle varies  monotonically between $0$ (for $t=1$, when we have a complete, finite-volume structure as in the end of Section \ref{sec: octahedron}) to $2\pi$ (for $t\to 0^+$). Analogously to the punctured torus of the previous section, these hyperbolic structures, when conjugate by $\gamma_t$ to real projective structures, have a limiting Euclidean structure for $t=0$, which is obtained by identifying opposite faces of the parallelepiped using translations. This can be proved by direct computations, very similar to those given in Section \ref{subsec:2d eucl}. In this Euclidean limit the structure extends to the link (i.e. there are no cone singularities on each link component, since the Euclidean angles equal $2\pi$), hence we recover a Euclidean three-torus.

The same holds on the spherical side: the spherical structures obtained by the same identifications on two copies of $\mathcal O(t)$, for $t<0$, converge to the Euclidean structure if rescaled by $\gamma_{|t|}$ in projective space. Cone singularities are larger than $2\pi$ for the spherical structures (as it can be again expected, for instance if one applies a Gauss-Bonnet formula to the spherical polygon obtained as the section of $\mathcal O(t)$ with a plane orthogonal to one of the edges). 

\begin{remark}
It is interesting to note that some spherical structure can be obtained also from the polytope $\mathcal O(-1)$, pictured in the right of Figure \ref{fig:regen oct eucl}, which is essentially a spherical wedge on a geodesic interval. One can glue opposite faces from a single copy of $\mathcal O(-1)$, and obtain a spherical metric on the 3-sphere with cone singularities along the Hopf link, having in this example the same angle (and the same length), of value $\pi/2$.
\end{remark}

\section{...and to half-pipe and Anti-de Sitter structures}\label{sec:HP}

The purpose of this section is to present the geometric transition in the pseudo-Riemannian setting which involves hyperbolic structures, Anti-de Sitter structures and the so-called \emph{half-pipe structures}, introduced in \cite{danciger}. We will then provide deformations of the hyperbolic structures of Sections \ref{sec:quadrilateral} and \ref{sec: octahedron} to exhibit concrete examples in dimensions 2 and 3.

\subsection{A pseudo-Riemannian geometric transition}

The \emph{Anti-de Sitter space} of dimension $n$ is defined, similarly to the sphere and the hyperbolic space, as 
$$\AdS^n=\{x\in\R^{n+1}\,:\,q_{2,n-1}(x)=-1\}~.$$
Endowed with the restriction of the standard bilinear form of signature $(2,n-1)$, $\AdS^n$ is a Lorentzian manifold of constant sectional curvature $-1$. It is not simply connected, being homeomorphic to $D^{n-1}\times S^1$. The isometry group $\Isom(\AdS^n)$ is identified to the group $\O(2,n+1)$, acts transitively on $\AdS^n$, and moreover every linear isometry between $T_{p}\AdS^n$ and $T_q\AdS^n$ extends to a global isometry of $\AdS^n$.

As for hyperbolic space, totally geodesic hyperplanes of $\AdS^n$ are the intersections of $\AdS^n$ with linear hyperplanes $P$ of $\R^{n+1}$ (which is always non-empty). They can be of three types: \emph{spacelike} when the restriction of the Lorentzian metric to $P\cap\AdS^n$ is positive definite (these are intrinsically isometric to $\Hyp^{n-1}$), \emph{timelike} when the restriction is still Lorentzian (intrinsically isometric to $\AdS^{n-1}$), or \emph{lightlike} when the restriction is a degenerate bilinear form. 

The most common model for Anti-de Sitter space is the projective model, namely $\AdS^n$ is the region defined by $q_{2,n-1}(x)<0$ in the projective sphere $\SP^n$. In the affine chart $x_0=1$, which however covers only a portion of $\AdS^n$, $\AdS^n$ is visualised as the interior of a one-sheeted hyperboloid. The projective model is also useful to see its \emph{boundary}, which is 
$$\partial\AdS^n=\{[x]\in\SP^n\,:\,q_{2,n-1}(x)=0\}~,$$
and is homeomorphic to $S^{n-2}\times S^1$.

We will again produce examples of geometric transitions in the sense of $(G,X)$--structures, starting from hyperbolic and landing on Anti-de Sitter structures, namely $(\Isom(\AdS^n),\AdS^n)$--structures. Let us see how this can happen, which will also motivate the definition of the intermediate geometry between $\Hyp^n$ and $\AdS^n$, called \emph{half-pipe}.

\begin{figure}[htb]
\includegraphics[height=7cm]{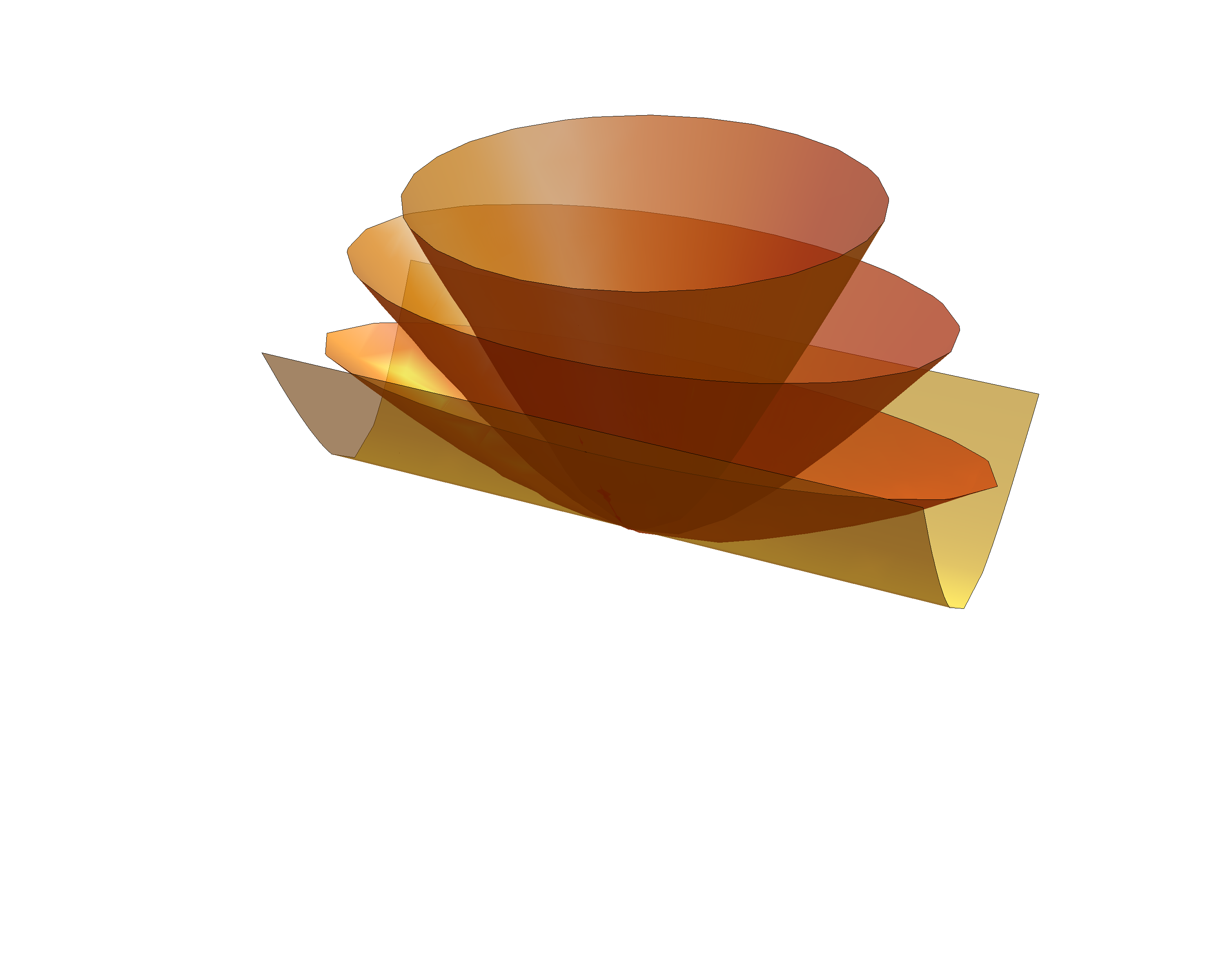}
\vspace{-2.6cm}
\caption{\footnotesize The Hausdorff limit of $\eta_t(\Hyp^n)$ in $\R^{n+1}$ is the so-called half-pipe space $\HP^n$ (in this picture $n=2$).}\label{fig:hausdorff2}
\end{figure}

By a slight inconsistency with \eqref{eq:quadratic forms}, we shall permute coordinates and denote by $q_{2,n-2}$ the quadratic form
$$q_{2,n-2}(x)=-x_0^2+x_1^2+\ldots+x_{n-1}^2-x_n^2$$
on $\R^{n+1}$. We also introduce for $t>0$ the linear transformation
\begin{equation} \label{eq etat}
\eta_t=\mathrm{diag}\left(1,\ldots,1,\frac{1}{t}\right)~,
\end{equation}
which will play in this section the role the transformations $\gamma_t$ had for the hyperbolic-Euclidean-spherical transition (Equation \eqref{eq gammat}). In fact, the Hausdorff limit of the closed hypersurfaces $\eta_t(\Hyp^n)$ in $\R^{n+1}$ is the so-called \emph{half-pipe space}
$$\HP^n=\{x=(x_0,\ldots,x_n)\,:\,-x_0^2+x_1^2+\ldots+x_{n-1}^2=-1,\,x_0>0\}~.$$
The last coordinate $x_n$ has been deliberately omitted in the definition of $\HP^n$, which is in fact the product of $\Hyp^{n-1}\subset \R^n$ and $\R$. Analogously to what happened in the spherical-to-Euclidean transition, the Hausdorff limit of $\eta_t(\AdS^n)$ as $t\to 0^+$ is $\pm\HP^n$. See Figure \ref{fig:hausdorff2}.

\vspace{-0.5cm}

If in the discussion of Section \ref{sec:eucl} we have observed that not only $\E^n$ is the Hausdorff limit of the rescaled copied of $\Hyp^n$ and $\Sph^n$, but also $\Isom(\E^n)$ is the Chabauty limit of the $\gamma_t$-conjugates of $\Isom(\Hyp^n)$ and $\Isom(\Sph^n)$, here one can \emph{compute} the limits of $\Isom(\Hyp^n)$ and $\Isom(\AdS^n)$, now conjugate by $\eta_t$, and \emph{define} their limits (which coincide) as the natural group acting on $\HP^n$.
This is how the half-pipe transformation group was first introduced in \cite{danciger}, namely the group of matrices of the form
\begin{equation} \label{eq:HPmatrices}
\begin{pmatrix} 
A &  0 \\
v & \pm 1 
\end{pmatrix}
\end{equation}
for $A\in\O_+(1,n-1)\cong\Isom(\Hyp^{n-1})$
and $v\in\R^n$. Indeed, it is easily checked that matrices of this form preserve the degenerate bilinear form $-x_0^2+x_1^2+\ldots+x_{n-1}^2$, preserve the sign of the $x_0$-coordinate, and therefore leave $\HP^n\subset\R^{n+1}$ invariant.

Of course half-pipe geometry have a projective model, since we can see $\HP^n$ as the subset of $\RP^n$ (or $\SP^n$) defined by $-x_0^2+x_1^2+\ldots+x_{n-1}^2<0$ (and $x_0>0$), and half-pipe transformations induce projective transformations in this model. Hence we can give the following definition of geometric transition in the setup of real projective structures (compare with Definition \ref{defi transition eucl}):

\begin{defi} \label{defi transition HP}
A \emph{geometric transition} on a manifold $\mathcal M$ from hyperbolic to Anti-de Sitter geometry, through half-pipe geometry, is a continuous path of real projective structures $\mathscr P_t$ on $M$, defined for $t\in (-\epsilon,\epsilon)$, which is conjugate to
\begin{itemize}
\item Hyperbolic structures for $t>0$;
\item Half-pipe structures for $t=0$;
\item Anti-de Sitter structures for $t<0$.
\end{itemize}
\end{defi}
Clearly by half-pipe structures we mean $(G,X)$--structures where $X=\HP^n$ and $G$ is the group of half-pipe transformations of the form \eqref{eq:HPmatrices}. In the next section we shall explore more deeply this geometry, and see that in fact we can think at $G$ as the group of isometries of Minkowski space.

\subsection{The geometry of half-pipe space and Minkowski space} \label{subsec:geometryHPM}
As a consequence of the above observation that the matrices of the half-pipe group \eqref{eq:HPmatrices} preserve the degenerate bilinear form $-x_0^2+x_1^2+\ldots+x_{n-1}^2$, we can endow $\HP^n$ with a natural \emph{degenerate} metric (by  restriction), which makes $\HP^n$ a pseudometric space. The degenerate directions are those tangent to the $\R$ factor in the decomposition $\HP^n=\Hyp^{n-1}\times\R$, that is, generated by $\partial/\partial{x_n}$.
Here of course \emph{natural} means a notion which is preserved by the group of transformations. 

There is however much more geometry we can introduce on $\HP^n$. Indeed, the matrices \eqref{eq:HPmatrices} have unit determinant, which permits to introduce a natural \emph{volume form}. A notion of \emph{length} is meaningful along the degenerate directions, since the bottom-right entry in \eqref{eq:HPmatrices} is also unit.
Finally, the notion of \emph{lines} makes sense by the projective model of $\HP^n$, which  also provides the \emph{boundary} of half-pipe space, namely
$$\partial\HP^n=\{[x]\in\RP^n\,:\,-x_0^2+x_1^2+\ldots+x_{n-1}^2=0\}~.$$

All these structures are interpreted in a transparent way in terms of Minkowski geometry. Let us denote by $\M^n$ the $n$-dimensional Minkowski space, namely $\R^n$ endowed with the Lorentzian metric $\langle\cdot,\cdot\rangle$ associated to $q_{1,n-1}$. Then we have a homeomorphism
\begin{equation}\label{eq:identificationHP}
\HP^n\cong\{\text{spacelike hyperplanes in }\M^n\}
\end{equation}
which associates to $x=(\bar x,x_n)\in \HP^n$ the hyperplane 
$$\mathcal F(x)=\{ y=(y_0,\ldots,y_{n-1})\in\M^n\,:\,\langle \bar x,y\rangle=x_n\}~.$$
In short, $\bar x=(x_0,\ldots,x_{n-1})\in\Hyp^n$ determines the future unit normal vector of $\mathcal F(x)$, while the last coordinate $x_n$ is the  signed distance of $\mathcal F(x)$ to the origin along the normal direction. Then we analogously have

\begin{equation}\label{eq:identificationHPbdy}
\partial\HP^n\cong\{\text{lightlike hyperplanes in }\M^n\}~.
\end{equation}

\begin{figure}[htb]
\includegraphics[height=6cm]{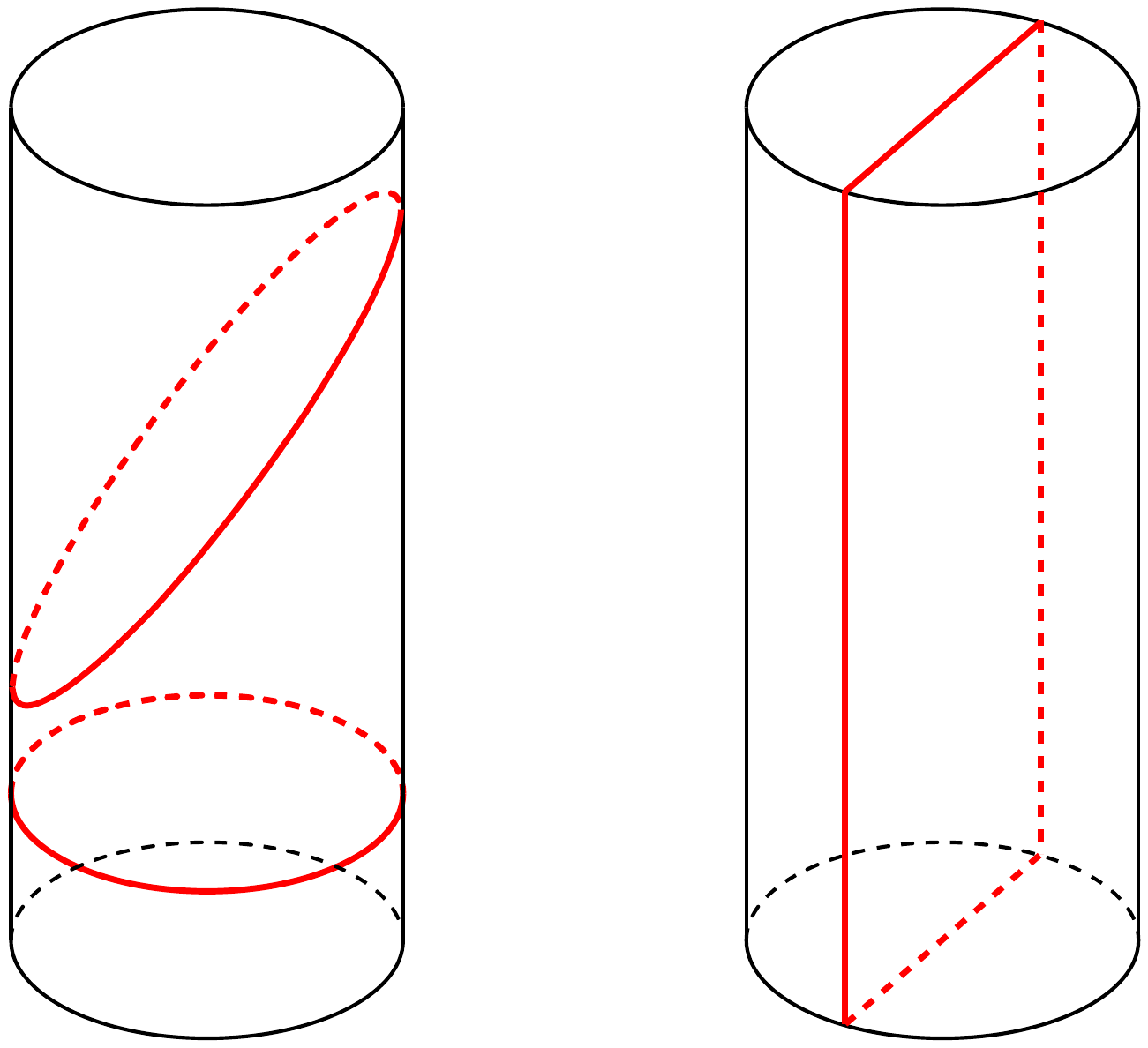}
\caption{\footnotesize Hyperplanes in $\HP^n$, which is an infinite cylinder in the projective model: on the left, of non-degenerate type, on the right, degenerate ones.}\label{fig:hyperplanesHP}
\end{figure}

\begin{remark} \label{rmk duality}
The homeomorphism $\mathcal F$ in \eqref{eq:identificationHP} can be thought as a \emph{duality} between $\HP^n$ and $\M^n$. In fact, we also have the dual correspondence
\begin{equation}\label{eq:identificationHPdual}
\M^n\cong\{\text{non-degenerate hyperplanes in }\HP^n\}~,
\end{equation}
where a \emph{non-degenerate hyperplane} is a half-pipe hyperplane, that is the intersection of $\HP^n$ with a linear hyperplane of $\R^{n+1}$, transverse to the degenerate direction. The map in \eqref{eq:identificationHPdual} associates to $v\in\M^n$ the set of all spacelike hyperplanes in $\M^n$ which contain $p$. These correspond indeed to a set of points in $\HP^n$, which turns out to be a non-degenerate hyperplane in $\HP^n$. Such non-degenerate hyperplanes in $\HP^n$ are isometric copies of $\Hyp^{n-1}$ for the half-pipe pseudometric. See Figure \ref{fig:hyperplanesHP}.
\end{remark}

It is easily checked (see \cite[Section 2.1.2]{barbotfillastre}, \cite{surveyseppifillastre}, \cite[Lemma 1.8]{rioloseppi}) that the homeomorphism \eqref{eq:identificationHP} induces a group isomorphism between $\Isom(\M^n)$ and the group of half-pipe transformations of the form \eqref{eq:HPmatrices}, by means of the obvious action of $\Isom(\M^n)$ on the space of spacelike hyperplanes of $\M^n$. The action extends to the boundary by \eqref{eq:identificationHPbdy}. Recalling that $\Isom(\M^n)\cong \O(1,n-1)\rtimes\R^n$, under this group isomorphism we have that:
\begin{itemize}
\item Linear isometries of $\M^n$ of the form $y\mapsto Ay$, for $A\in \O_+(1,n-1)$ (namely, $A$ preserves the sign of the $y_0$--coordinate) correspond to matrices $\begin{pmatrix} 
A &  0 \\
0 & 1 
\end{pmatrix}$;
\item The linear isometry $y\mapsto -y$, which generates the center of $\O(1,n-1)$ (in fact one has $\O(1,n-1)\cong\O_+(1,n-1)\times(\Z/2\Z)$), corresponds to 
$\begin{pmatrix} 
\mathrm{id} &  0 \\
0 & -1 
\end{pmatrix}$;
\item Translations $y\mapsto y+v$ correspond to 
$\begin{pmatrix} 
\mathrm{id} &  0 \\
v^T\!J & 1 
\end{pmatrix}$, where $J=\mathrm{diag}(-1,1,\ldots,1)$.
\end{itemize}
Three important concrete examples are the following. See also Figure \ref{fig:hyperplanesHP}.
\begin{enumerate}
\item First, an isometry of the form $y\mapsto -y+v$ for $v\in\M^n$ fixes the point $v/2\in\M^n$ and leaves invariant all the spacelike hyperplanes which contain $v/2$. This means that the corresponding action on $\HP^n$ fixes (pointwise) a {non-degenerate hyperplane} (see Remark \ref{rmk duality}). As the Minkowski isometry $y\mapsto -y+v$  reverses the degenerate direction in $\HP^n$, it induces a transformation of $\HP^n$ which plays the role of a \emph{reflection in a non-degenerate hyperplane}.
\item Second, two non-degenerate hyperplanes $\mathcal P$ and $\mathcal P'$ in $\HP^n$, which correspond to two points $v$ and $v'$ in $\M^n$ (as in \eqref{eq:identificationHPdual}), intersect if and only if $v-v'$ is spacelike. In this case, the composition of the two reflections along $\mathcal P$ and $\mathcal P'$, as in the previous point, give the analogue of a \emph{rotation} in $\HP^n$, which fixes pointwise $\mathcal P\cap\mathcal P'$ (a copy of $\Hyp^{n-2}$). In $\Isom(\M^n)$, such rotation corresponds to a translation in the direction of $v-v'$.
\item Finally, {degenerate hyperplanes} in $\HP^n$, that is hyperplanes which contain a degenerate direction, correspond to all spacelike hyperplanes $\mathcal P$ in $\M^n$ whose normal vector lies in a given timelike hyperplane $\mathcal T$ (hence in a geodesic of $\Hyp^{n-1}$). The isometries of $\M^n$ whose linear part is a reflection in $\mathcal T$, and whose translation part  is orthogonal to $\mathcal T$, leave setwise invariant each of these spacelike hyperplanes $\mathcal P$. These all act on $\HP^n$ as \emph{reflections in the degenerate hyperplane}. In particular, half-pipe reflections in degenerate hyperplanes are not unique, but they constitute one-parameter groups.
\end{enumerate}

\subsection{On a (punctured) torus} \label{subsec:2d HP}
Let us now go back to the punctured torus of Section \ref{sec:quadrilateral}, and produce a different deformation which gives rise to a geometric transition to half-pipe and Anti-de Sitter geometry. This example already appears in a similar form in \cite[Section 3.3]{dancigertransition}. The idea is to deform the quadrilateral by keeping the length of one axis fixed, and letting the other axis collapse to zero. 

\begin{figure}[htb]
\includegraphics[width=\textwidth]{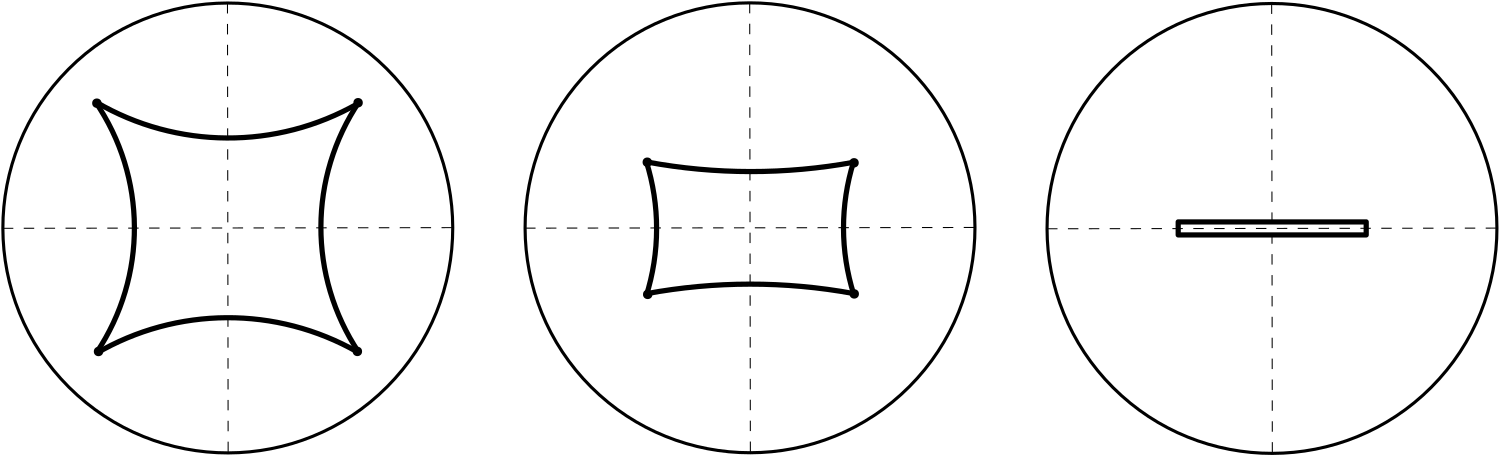}
\caption[The cuboctahedron.]{\footnotesize The family of quadrilaterals $\mathcal Q'(t)$, collapsing to a geodesic segment.}\label{fig:quadrilateral3}
\end{figure}

Explicitly, consider the quadrilateral $\mathcal Q'(t)\subset\Hyp^2$ for $t>0$, which is defined as the intersection of the four half-planes $(-1:\pm \sqrt 2:0)$ and $(-t:0:\pm\sqrt{1+t^2})$, seen as elements of $\SP^2$ as usual. See Figure \ref{fig:quadrilateral3}. When $t=1$, this is exactly the ideal quadrilateral described in Section \ref{sec:quadrilateral}, with vertices \eqref{eq:vertices ideal quad}, and when $t\in(0,1)$ a quadrilateral with vertices in $\Hyp^2$, where the distance between the top and bottom sides is $2\mathrm{arcsinh}(t)$. Also, a simple computation shows that the angle at the vertices of $\mathcal Q'(t)$ takes the value $\mathrm{arccos}(t)$, hence tends to $\pi/2$ as $t\to 0^+$.
We shall also extend the definition of the quadrilateral $\mathcal Q'(t)$ for negative $t$ to an Anti-de Sitter quadrilateral, in which case the two ``vertical'' sides are timelike. We will thus declare that for $t\in (-1,0)$ the quadrilateral $\mathcal Q'(t)$ is the intersection of the four half-planes $(-1:\pm \sqrt 2:0)$ and $(-|t|:0:\pm\sqrt{1-t^2})$, which makes the timelike distance between the two ``horizontal''  geodesics (the non-constant ones, which are spacelike) equal to $2\mathrm{arcsin}(|t|)$. 

\begin{figure}[htb]
\centering
\begin{minipage}[c]{.25\textwidth}
\centering
\includegraphics[height=4.5cm]{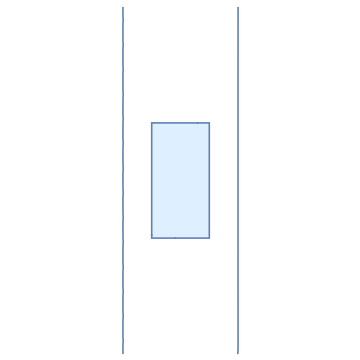}
\end{minipage}%
\begin{minipage}[c]{.35\textwidth}
\centering
\includegraphics[height=4.5cm]{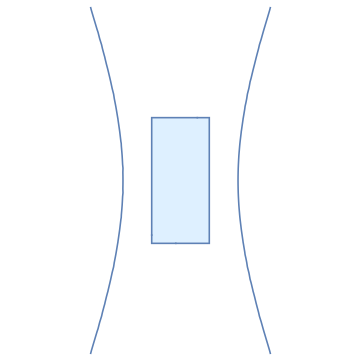}
\end{minipage}%
\begin{minipage}[c]{.3\textwidth}
\centering
\includegraphics[height=4.5cm]{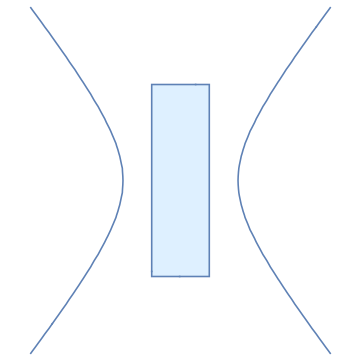}
\end{minipage}
\caption{\footnotesize The continuation of $\mathcal Q'(t)$, rescaled by $\eta_{|t|}$, from a half-pipe quadrilateral to the Anti-de Sitter side. When $t\to -1$, the vertical sides become lightlike and the quadrilateral gets out of an affine chart.}
\label{fig:regen quad HP}
\end{figure}

Rescaling in the direction $x_2$ by means of the transformations $\eta_t$ (Equation \eqref{eq etat}), we obtain a family of quadrilaterals $\eta_t(\mathcal Q'(t))$ in projective space which converges as $t\to 0^\pm$ to a half-pipe quadrilateral whose sides are defined by the following elements of $\SP^2$: $(-1:\pm \sqrt 2:0)$, which are degenerate for the geometry of $\HP^2$, and $(-1:0:\pm 1)$, which are non-degenerate. See Figure \ref{fig:regen quad HP}. In fact, the half-planes $(-1:\pm \sqrt 2:0)$ of $\mathcal Q'(t)$ are clearly left invariant by $\eta_t$, while in general the effect of $\eta_t$ on a linear form $\alpha$ amounts to multiplying all the coordinates of the corresponding $\alpha$ by $1/t$, except the last one, which makes evident that the rescaled limit of the other two half-spaces is the one we claimed.

One still produces hyperbolic structures $\mathscr H'(t)$  on the punctured torus, glueing the sides two by two using the compositions of hyperbolic reflections in the axes of symmetry and in the sides themselves, and analogously Anti-de Sitter structures $\mathscr A'(t)$ for negative $t$. These structures collapse to a circle as $t\to 0^\pm$, and we claim that they provide a geometric transition as in Definition \ref{defi transition HP}. We shall show that if $\mathscr P'_t$ are the real projective structures obtained by applying $\eta_t$ to $\mathscr H'(t)$, and $\eta_{|t|}$ to $\mathscr A'(t)$ when $t$ is negative, then $\mathscr P'_t$ converges to the same half-pipe structure when $t\to 0^\pm$.

Clearly one can produce a family of charts for $\mathscr P'_t$ which vary continuously in $t$, all of which have image with non-empty intersection with $\eta_{|t|}(\mathcal Q'(t))$, similarly to Figure \ref{fig:quadrilateral2}. The important part is to show that the change of coordinates pass to the limit, which is less evident then in the case of geometric transition to Euclidean geometry. In fact, while in hyperbolic  and Euclidean space there is a unique reflection in a geodesic (or a hyperplane, in general), we have seen that the same does not hold in half-pipe geometry (see Section \ref{subsec:geometryHPM}, last item). For this reason, in the limit half-pipe quadrilateral we shall somehow \emph{choose} a reflection in the ``vertical'' sides (that are degenerate), which is not simply determined by the quadrilateral itself.

\begin{figure}[htb]
\includegraphics[height=4.5cm]{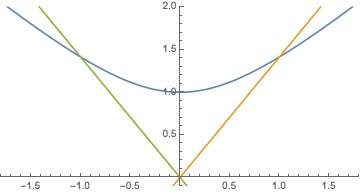}
\caption{\footnotesize The boost $\mathcal R$ in $\M^2$ maps the left line to the right line, and translates by a length of $2\mathrm{arcsinh}(1)$ along the hyperbola $\Hyp^1$.}\label{fig:minkowski11}
\end{figure}

To exhibit the identifications, we will of course study how the hyperbolic reflections pass to the limit when conjugate by $\eta_t$. It is easy to check that the reflections in the ``vertical'' sides of $\mathcal Q(t)$, which do not depend on $t$ and in fact have the same expression for hyperbolic and Anti-de Sitter space, commute with $\eta_t$, and so does the reflection in the vertical axis of symmetry. In the limit, we obtain half-pipe transformations  which correspond under \eqref{eq:identificationHP} to linear isometries of $\M^2$ which are Minkowski reflections in the timelike lines $y_0=0$ and $y_0=\pm \sqrt 2 y_1$, using orthonormal coordinates $(y_0,y_1)$  where $\partial_{y_0}$ is a negative direction and $\partial_{y_1}$ a positive one, consistently with the notation of Section \ref{subsec:geometryHPM}. This falls indeed in item (3) of that section.  Hence the identification between the left and right side corresponds to a Minkowski boost $\mathcal R$, which fixes the origin and translates along the hyperbola $\Hyp^1\subset\M^2$ by an intrinsic length of $2\mathrm{arcsinh}(1)$.  See Figure \ref{fig:minkowski11}.

Next, let us consider the identification between top and bottom sides. Like before, the reflection in the horizontal axis of symmetry is independent of $t$, has the common expression $(x_0,x_1,x_2)\mapsto(x_0,x_1,-x_2)$ for each of the geometries $\Hyp^2$, $\HP^2$ and $\AdS^2$, and in fact commutes with $\eta_t$. Under the isomorphism with $\Isom(\M^2)$, the half-pipe reflection corresponds to $y\mapsto -y$, which is a reflection as in item (1) of Section \ref{subsec:geometryHPM}, with trivial translation part. We are only left with understanding the rescaled limit of the hyperbolic and Anti-de Sitter reflections in the vertical sides, which do depend on $t$. Using the expression of Equation \eqref{eq:hyp reflection} and repeating a similar computation with $\eta_t$ instead of $\gamma_t$, we find that the limit in this case is the transformation $$\begin{pmatrix} 
\mathrm{id} &  0 \\
v & -1 
\end{pmatrix}$$
for $v=(\pm 1,0)$. In $\Isom(\M^2)$, this corresponds to the transformation $y\mapsto -y+v$ as expected from item (1) of Section \ref{subsec:geometryHPM}. Thus the identification between top and bottom side corresponds to a Minkowski translation of vector $(1,0)$.

As a final observation, the half-pipe structure we constructed on the punctured torus do not extend to a honest structure on the torus, differently to what happened in the hyperbolic-to-Euclidean transition. In fact, one can compute the holonomy around the puncture (which is exactly the obstruction to extending the geometric structure) as the commutator of the half-pipe transformations which are used to identify the sides of the quadrilateral. In $\Isom(\M^2)$, such a commutator is of the form $y\mapsto y+\mathcal Rv-v$, where $\mathcal R$ is a boost as above. Since the vector $v=(1,0)$ is timelike, $\mathcal Rv-v$ is spacelike, and we have a \emph{half-pipe rotation} as described in item (2) of Section \ref{subsec:geometryHPM}: the punctured torus has the analogue of a \emph{cone singularity} in half-pipe geometry.

\begin{figure}[htb]
\includegraphics[height=4cm]{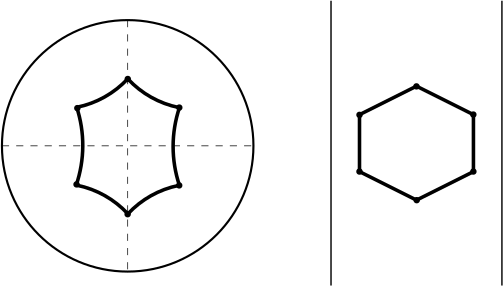}
\caption{\footnotesize A schematic picture of the hyperbolic hexagon and the corresponding half-pipe limit hexagon, obtained by ``switching'' left and right side of the quadrilaterals.}\label{fig:hexagon}
\end{figure}

This fact, which is a manifestation of the general fact that there are no non-singular half-pipe structures on compact manifolds (\cite[Section 3.4]{danciger}), might seem counterintuitive but is more easily visualised if one does a cut-and-paste on the quadrilateral $\mathcal Q(t)$, switching the left and right half as in Figure \ref{fig:hexagon}. Now the punctured torus is obtained from a hyperbolic (or Anti-de Sitter) hexagon, which converges to a segment as $t\to 0^\pm$, and to a half-pipe hexagon if rescaled by $\eta_{|t|}$.

\subsection{On the Borromean rings complement and the three-torus} \label{subsec:3d HP}
The purpose of this section is to describe a family of three-dimensional polytopes $\mathcal O'(t)$ which deform the octahedron of Figure \ref{fig:octahedron} (hence $\mathcal O'(1)=\mathcal O$) and collapse, differently from Section \ref{subsec:3d eucl}, to a polygon in a totally geodesic plane of $\Hyp^3$. The family can actually be continued for negative times to $\AdS^3$ and, after rescaling, will transition through a half-pipe polytope. Similarly to Section \ref{subsec:3d eucl}, the important feature of $\mathcal O'(t)$ is that, in some sense, many of the dihedral angles are right.

To construct such deformation, we will impose that $\mathcal O'(t)$ maintains the symmetries by reflection of the right-angles ideal octahedron $\mathcal O$. These, we recall, are expressed by changing sign to $x_1$, $x_2$ or $x_3$, and in fact are isometries not only for $\Hyp^3\subset\R^4$ but also for $\AdS^3$ and $\HP^3$. We will however break the symmetry by rotations of angle $\pi/2$ in a vertical axis. Moreover, we will require that the intersection of $\mathcal O'(t)$ with the totally geodesic plane $\{x_3=0\}$ is constant in $t$, and is therefore the regular ideal quadrilateral as for $\mathcal O$.

\begin{table}[htb]
\begin{eqnarray*}
\left( -|t|:-\sqrt{2}|t|:0:-1 \right) , & &
\left( -1:-\sqrt{2}:0:+t \right),\\
\left( -|t|:0:-\sqrt{2}|t|:+1 \right),& &
\left( -1:0:-\sqrt{2}:-t \right),\\
\left( -|t|:+\sqrt{2}|t|:0:-1 \right),& &
\left( -1:+\sqrt{2}:0:+t \right),\\
\left( -|t|:0:+\sqrt{2}|t|:+1 \right),& &
\left( -1:0:+\sqrt{2}:-t \right).\\
\end{eqnarray*}
\caption{\footnotesize The half-spaces defining the deformation $\mathcal O'(t)$, which is hyperbolic for $t>0$ and Anti-de Sitter for $t<0$. When $t=0$, the half-spaces of the left column become $\pm x_3>0$, while those of the right column become orthogonal to the plane $x_3=0$ and in fact coincide which bound $\mathcal Q$ (adding a zero in the last entry). Hence $\mathcal O'(0)$ degenerates to a regular ideal quadrilateral. For negative $t$, the boundary of the half-spaces in the left column is spacelike, and timelike in the right column.}\label{table:link}
\end{table}

The list of half-spaces defining $\mathcal O'(t)$ is given in Table \ref{table:link}, for $t\in (-1,1)$. For $t=1$, this coincides with Table \ref{table:oct ideal}. The reader might want to compare also with Table \ref{table:oct collapse1} to see the differences in the two deformations. Each of the eight faces of $\mathcal O'(t)$ is contained in the boundary of one of the eight defining half-spaces in Table \ref{table:link}, and they come in pairs, which appear in the same row. The property of each pair is that the corresponding faces are orthogonal, meet in a geodesic of the ideal quadrilateral $\mathcal O'(t)\cap\{x_3=0\}$, and lie one above and one below the totally geodesic plane $\{x_3=0\}$. Orthogonality holds for the hyperbolic metric when $t>0$, and for the Anti-de Sitter metric for $t<0$, which is easily checked since the two metrics differ by the sign in the last coordinate, and in fact when $t$ becomes negative there is a change of sign only in the last entry in the right column. 

\begin{figure}[htb]
\centering
\begin{minipage}[c]{.25\textwidth}
\centering
\includegraphics[scale=0.25]{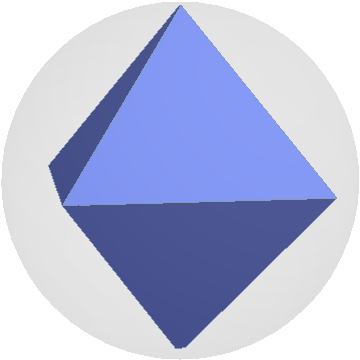}
\end{minipage}%
\begin{minipage}[c]{.25\textwidth}
\centering
\includegraphics[scale=0.25]{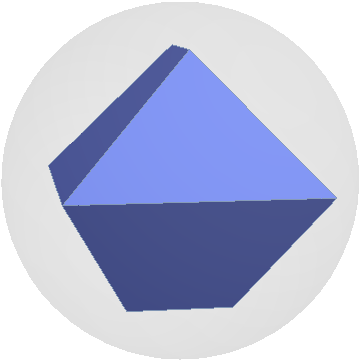}
\end{minipage}%
\begin{minipage}[c]{.25\textwidth}
\centering
\includegraphics[scale=0.25]{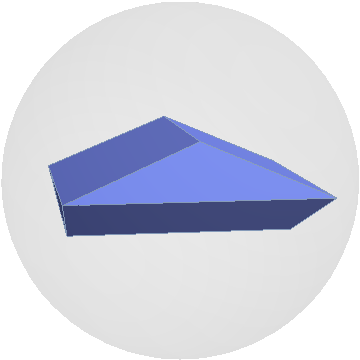}
\end{minipage}
\begin{minipage}[c]{.24\textwidth}
\centering
\includegraphics[scale=0.25]{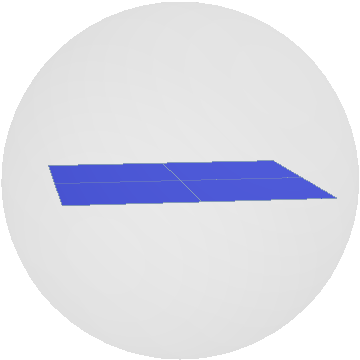}
\end{minipage}%
\caption{\footnotesize The deformation ${\mathcal O}(t)$ in an affine chart (Klein model of $\Hyp^3$), from the right-angled ideal octahedron to the collapse at a regular ideal quadrilateral.}
\label{fig:collapse oct HP}
\end{figure}

\begin{figure}[htb]
\centering
\begin{minipage}[c]{.33\textwidth}
\centering
\includegraphics[scale=0.3]{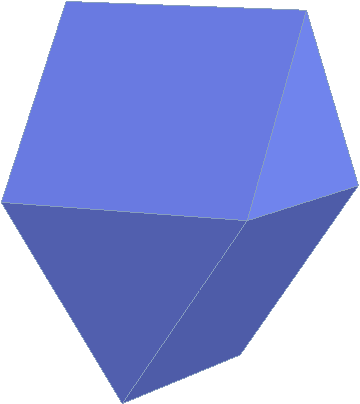}
\end{minipage}%
\begin{minipage}[c]{.33\textwidth}
\centering
\includegraphics[scale=0.28]{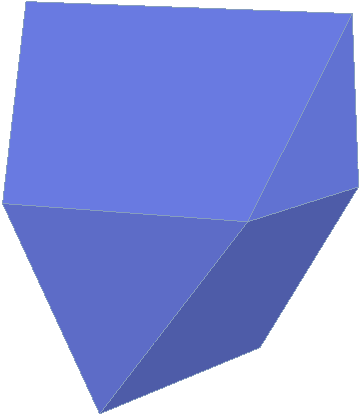}
\end{minipage}%
\begin{minipage}[c]{.33\textwidth}
\centering
\includegraphics[scale=0.33]{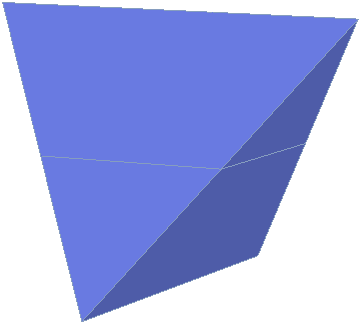}
\end{minipage}
\caption{\footnotesize The rescaled polytope $\eta_{|t|}\mathcal O'(t)$ for negative $t$, which interpolates from a half-pipe to Anti-de Sitter geometry. When $t=-1$ we have a polytope in $\AdS^3$ with four  lightlike faces.}
\label{fig:regen oct HP}
\end{figure}

There are new edges which appear in $\mathcal O'(t)$ when $t<1$, and these are the only non-right dihedral angles. Each face in the upper half-space of $\Hyp^3$ (that is, with negative last coordinate) does indeed intersect orthogonally two of the other faces in the upper half-space, and form a dihedral angle with the fourth one. Such dihedral angle varies monotonically between $0$ and $\pi$. See Figure \ref{fig:collapse oct HP}. Of course the same holds symmetrically on the lower half-space, and on the Anti-de Sitter side of the deformation, with the additional remark that the planes bounding the half-spaces of Table \ref{table:link} are spacelike on the left column and timelike on the right column. When $t=-1$, the polyhedron degenerates to a ``wedge'' with lightlike faces, see Figure \ref{fig:regen oct HP} on the right.

When $t\to 0^\pm$, the qualitative behaviour of the half-spaces in the left and in the right column is very different. In fact, the former converge to $(0:0:0:\pm 1)$, namely to the half-space $x_3\geq 0$ or $x_3\leq 0$, while the latter converge to half-spaces orthogonal to the plane $x_3=0$ and have intersection with it a half-plane which bounds the ideal quadrilateral $\mathcal Q$. Such difference is highlighted when one considers the rescaled limit of $\eta_{|t|}\mathcal O'(t)$ as $t\to 0^\pm$. To compute the rescaled limit, observe that the action of $\eta_{|t|}$ on $\SP^3$ induces on $\SP^{3,*}$ an action given by multiplying the last coordinate by $|t|$, or equivalently by multiplying all the coordinates by $1/|t|$ except the last one. Hence $\eta_{|t|}\mathcal O'(t)$ converges to a half-pipe polytope, which is the intersection of the half-spaces in Table \ref{table:link rescaled} and is pictured in Figure \ref{fig:regen oct HP} on the left. Again we have two types of planes: four degenerate (the ``vertical'' ones) and four non-degenerate.

\begin{table}[htb]
\begin{eqnarray*}
\left( -1:-\sqrt{2}:0:-1 \right) , & &
\left( -1:-\sqrt{2}:0:0 \right),\\
\left( -1:0:-\sqrt{2}:+1 \right),& &
\left( -1:0:-\sqrt{2}:0 \right),\\
\left( -1:+\sqrt{2}:0:-1 \right),& &
\left( -1:+\sqrt{2}:0:0 \right),\\
\left( -1:0:+\sqrt{2}:+1 \right),& &
\left( -1:0:+\sqrt{2}:0 \right).\\
\end{eqnarray*}
\caption{\footnotesize The half-spaces defining the limiting half-pipe polytope. The four faces corresponding to the left column are degenerate, those of the right column are non-degenerate.}\label{table:link rescaled}
\end{table}

We can construct geometric structures (hyperbolic and Anti-de Sitter) starting from the polytopes $\mathcal O'(t)$, for instance on the Borromean rings complement. 
In complete analogy with the previous examples, to show that these geometric structures transition through a half-pipe structure according to Definition \ref{defi transition eucl}, we have to show that reflections in the faces of $\mathcal O(t)$, as elements of $\Isom(\Hyp^3)$ and $\Isom(\AdS^3)$, converge to half-pipe reflections when conjugate by $\eta_{|t|}$. To see this, we shall distinguish two cases, namely the half-spaces in the left and right column. A direct computation, using the isomorphism \eqref{eq:identificationHP} and in the same spirit as Section \ref{subsec:2d HP}, show that the limiting reflections correspond to isometries in $\Isom(\M^3)$ such that:
\begin{itemize}
\item The linear part is $-\mathrm{id}$, for the four elements in the left column;
\item For the four elements in the right column, the linear part is a reflection in the corresponding geodesic side of the ideal quadrilateral $\mathcal Q\subset\Hyp^2$;
\item The translation part is $\pm(1,-\sqrt 2,0)$ for the elements in the first row, where the sign depends on the sign of the last coordinate in Table \ref{table:link}, and similarly for the following rows.
\end{itemize}
Half-pipe geometric structures with cone singularities are then constructed, so as to produce a geometric transition to half-pipe geometry as in Definition \ref{defi transition eucl}. For instance, one can verify that the reflections in sides of the half-pipe polytope which come from right dihedral angles of $\mathcal O'(t)$ commute, and therefore the geometric structure is non-singular along the corresponding edges, while in correspondence of the non-right dihedral angles one has cone singularities in the half-pipe sense (i.e. the holonomy is a ``half-pipe rotation'', see Section \ref{subsec:geometryHPM}).

\section{A glimpse into dimension four}\label{sec:4d}

In this final section we provide a quick overview on the four-dimensional construction of \cite{rioloseppi}, which gives examples of geometric transition from hyperbolic to half-pipe and Anti-de Sitter structures in dimension four. The starting point, similarly to the previous examples, is a family $\mathcal P(t)$ of deforming convex hyperbolic polytopes, which had been introduced by Kerckhoff and Storm in \cite{KS}. The list of defining vectors is given in Table \ref{table:walls}, where the parameter $t$ can also be taken to be negative (in $(-1,0)$) and in that case we will have an Anti-de Sitter polytope, and the extension to the Anti-de Sitter side was given in \cite{rioloseppi}. There are 22 defining half-spaces and they are collected in three groups: 8 of the form $\p i$, 8 of the form $\m i$, and 6 of the form $\l X$. In fact, Kerckhoff and Storm first obtained this deformation by removing two walls (two additional ``letters'' in the list) from the hyperbolic 24-cell, and observing, with the aid of a computer, that the 22 remaining hyperplanes admit deformations which preserve ``many'' orthogonality relations. Let us explain the qualitative behaviour of this polytope $\mathcal P(t)$.

\begin{table} [htb]
\begin{eqnarray*}
\p{0} = \left( -\sqrt{2}\ |t|:+|t|:+|t|:+|t|:+1 \right) , & &
\m{0} = \left( -\sqrt{2}:+1:+1:+1:-t \right),\\
\p{1} = \left( -\sqrt{2}\ |t|:+|t|:-|t|:+|t|:-1 \right),& &
\m{1} = \left( -\sqrt{2}:+1:-1:+1:+t \right),\\
\p{2} = \left( -\sqrt{2}\ |t|:+|t|:-|t|:-|t|:+1 \right),& &
\m{2} = \left( -\sqrt{2}:+1:-1:-1:-t \right),\\
\p{3} = \left( -\sqrt{2}\ |t|:+|t|:+|t|:-|t|:-1 \right),& &
\m{3} = \left( -\sqrt{2}:+1:+1:-1:+t \right),\\
\p{4} = \left( -\sqrt{2}\ |t|:-|t|:+|t|:-|t|:+1 \right),& &
\m{4} = \left( -\sqrt{2}:-1:+1:-1:-t \right),\\
\p{5} = \left( -\sqrt{2}\ |t|:-|t|:+|t|:+|t|:-1 \right),& &
\m{5} = \left( -\sqrt{2}:-1:+1:+1:+t \right),\\
\p{6} = \left( -\sqrt{2}\ |t|:-|t|:-|t|:+|t|:+1 \right),& &
\m{6} = \left( -\sqrt{2}:-1:-1:+1:-t \right),\\
\p{7} = \left( -\sqrt{2}\ |t|:-|t|:-|t|:-|t|:-1 \right),& &
\m{7} = \left( -\sqrt{2}:-1:-1:-1:+t \right),\\
\l{A} = \left( -1:+\sqrt{2}:0:0:0 \right),& &
\l{B} = \left( -1:0:+\sqrt{2}:0:0 \right),\\
\l{C} = \left( -1:0:0:+\sqrt{2}:0 \right),& &
\l{D} = \left( -1:0:0:-\sqrt{2}:0 \right),\\
\l{E} = \left( -1:0:-\sqrt{2}:0:0 \right),& &
\l{F} = \left( -1:-\sqrt{2}:0:0:0 \right).
\end{eqnarray*}
\caption{\footnotesize The half-spaces in $\SP^4$ that define the projective polytope ${\mathcal P}_t$, given as usual by elements of $\SP^{4,*}$. They are collected in three groups with different behaviour.}\label{table:walls}
\end{table}

\begin{figure}[htb]
\includegraphics[scale=.7]{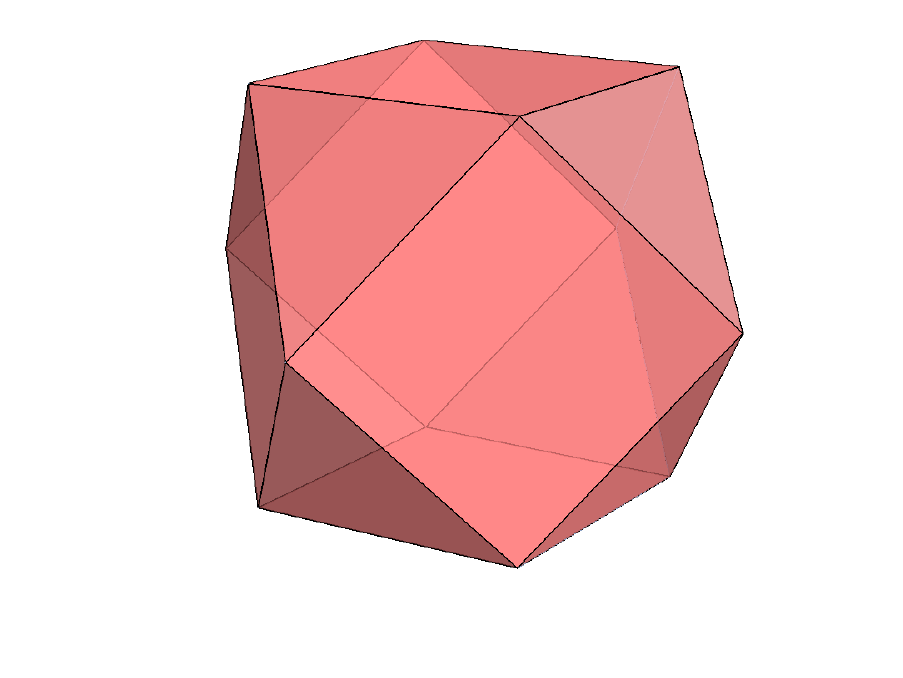}
\vspace{-1cm}
\caption{\footnotesize 
{The ideal right-angled cuboctahedron in $\Hyp^3$.}}\label{fig:cuboct}
\end{figure}

First of all, the intersection of $\mathcal P(t)$ with the ``horizontal'' totally geodesic hyperplane defined by $x_4=0$ (which is a copy of $\Hyp^3$, both in $\Hyp^4$ and in $\AdS^4$) is independent of $t$, and consists of the \emph{right-angled ideal cuboctahedron}, pictured in Figure \ref{fig:cuboct}.

The hyperplanes of the form $\p i$ are easily seen to converge to the horizontal hyperplane $x_4=0$ as $t\to 0^\pm$, and the corresponding half-spaces either to $x_4\leq 0$ or to $x_4\geq 0$. There are thus 8 hyperplanes with this behaviour, and the corresponding walls of $\mathcal P(t)$ lie in the upper half-space $x_4\geq 0$ for 4 of them, and in the lower half-space for the other four. When $t$ is negative, these walls are spacelike in $\AdS^4$.

The hyperplanes of the form $\m i$ instead converge to a hyperplane which is orthogonal to $x_4=0$ when $t$ approaches $0$. In fact, for $t$ negative these hyperplanes are timelike. Each wall $\m i$ intersects orthogonally, among others, the corresponding wall $\p i$, and their intersection lives in the hyperplane $x_4=0$ and is one of the triangular faces of the cuboctahedron (Figure \ref{fig:cuboct}).

Finally, the ``letter'' hyperplanes are constant (i.e. do not depend on $t$), and are timelike in $\AdS^4$ for negative times. They intersect the hyperplane $x_4=0$ in a quadrilateral face of the cuboctahedron. In fact, we can look at a ``slice'' of $\mathcal P(t)$, namely the intersection of $\mathcal P(t)$ with a hyperplane of the form $\l X$, which is a copy of $\Hyp^3$ (for $t>0$) or $\AdS^3$ (for $t<0$). Such intersection is a deforming three-dimensional polytope and is isometric exactly to the deformation $\mathcal O'(t)$ which we studied in Section \ref{subsec:3d HP} (see also Figures \ref{fig:collapse oct HP} and \ref{fig:regen oct HP}). The intersection with the horizontal plane is in fact a regular ideal quadrilateral, which is a face of the cuboctahedron, and in Figure \ref{fig:collapse oct HP} one sees that the walls of the form $\l X$ intersect 8 other walls of $\mathcal P(t)$: four of type $\p i$ (which become ``horizontal'' as $t\to 0^\pm$) and four of type $\m i$ (which become ``vertical'').  

\begin{figure}[htb]
\centering
\begin{minipage}[c]{.25\textwidth}
\centering
\includegraphics[scale=0.2]{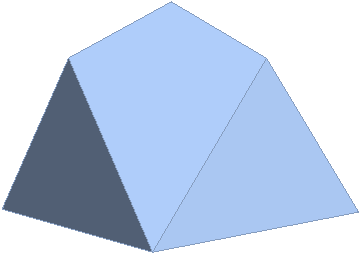}
\end{minipage}%
\begin{minipage}[c]{.25\textwidth}
\centering
\vspace{0.3cm}
\includegraphics[scale=0.2]{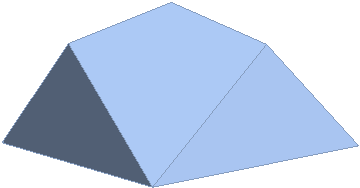}
\end{minipage}%
\begin{minipage}[c]{.25\textwidth}
\vspace{0.7cm}
\centering
\includegraphics[scale=0.2]{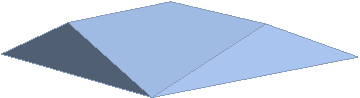}
\end{minipage}
\begin{minipage}[c]{.2\textwidth}
\vspace{1cm}
\centering
\includegraphics[scale=0.2]{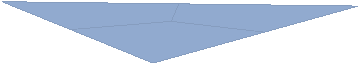}
\end{minipage}%
\caption{\footnotesize The intersection of $\mathcal P(t)$ with a hyperplane of the form $\m i$ is a hyperbolic/Anti-de Sitter polytope which collapses to an ideal triangle.}
\label{fig:collapse theta}
\end{figure}

Also the intersection of $\mathcal P(t)$ with a hyperplane of the form $\m i$ can be visualised, and is pictured in Figure \ref{fig:collapse theta} (in fact only the part lying in the half-space $x_4\geq 0$ is represented). Here one sees the intersections of $\m i$ with three walls of the form $\l X$ (the vertical ones, independent of $t$) and three of the form $\p i$ (becoming horizontal in the collapse), mutually orthogonal.

In \cite{rioloseppi} it is proved that, when rescaled by $\eta_{|t|}$, the polytope $\mathcal P(t)$ converges to a half-pipe polytope, for which the walls of the form $\l X$ and $\m i$ are degenerate, while those of the form $\p i$ are non-degenerate. This is used to produce a geometric transition hyperbolic/half-pipe/Anti-de Sitter on four-dimensional manifolds. More precisely, one starts from a complete, finite-volume hyperbolic three-manifold $\mathcal N$ which is obtained by glueing several copies of the ideal right-angled cuboctahedron of Figure \ref{fig:cuboct}. There are many ways to produce such $\mathcal N$, for instance one can first double the cuboctahedron along all quadrilateral faces, and then double again along all triangular faces, thus producing a cusped hyperbolic manifold out of four copies.

The key property of $\mathcal P(t)$ is that it is highly symmetric, in the sense that each self-isometry of the cuboctahedron can be extended to a self-isometry of $\mathcal P(t)$. (There is also a self-isometry which switches ``above'' and ``below'' with respect to the hyperplane $x_4=0$, which is however not a simple reflection, exactly like it happened for $\mathcal O'(t)$ in Figure \ref{fig:collapse oct HP}.) Hence one can ``follow'' the glueings which are used to produce $\mathcal N$, so as to glue several copies of $\mathcal P(t)$ and produce a four-manifold $\mathcal M$, homeomorphic to $\mathcal N\times S^1$, and endowed with a geometric transition to half-pipe and Anti-de Sitter geometry. The resulting manifold $\mathcal M$ is not compact, but its ends are \emph{cusps}, a notion which can be extended naturally to half-pipe and Anti-de Sitter geometry. In particular, the hyperbolic and Anti-de Sitter metrics constructed in the transition have finite volume. For completeness, we mention that  the deformation is to be considered in the interval $(-1,1/\sqrt 3)$, since on the hyperbolic side there are changes in the combinatorics of $\mathcal P(t)$ when $t>1/\sqrt 3$, as studied in \cite{MR}.

Two final remarks are to be made here. First, to exhibit the limiting half-pipe structure it is not enough to start from the limiting half-pipe polytope, but it is also necessary to determine the reflections in its faces which have to be used in the glueings. Recall indeed from Section \ref{subsec:geometryHPM} that half-pipe reflections in degenerate hyperplanes are not unique. This is done exactly in the same spirit as in Sections \ref{subsec:2d HP} and \ref{subsec:3d HP}, although computationally more intricate. 

\begin{figure}[htb]
\includegraphics[scale=.18]{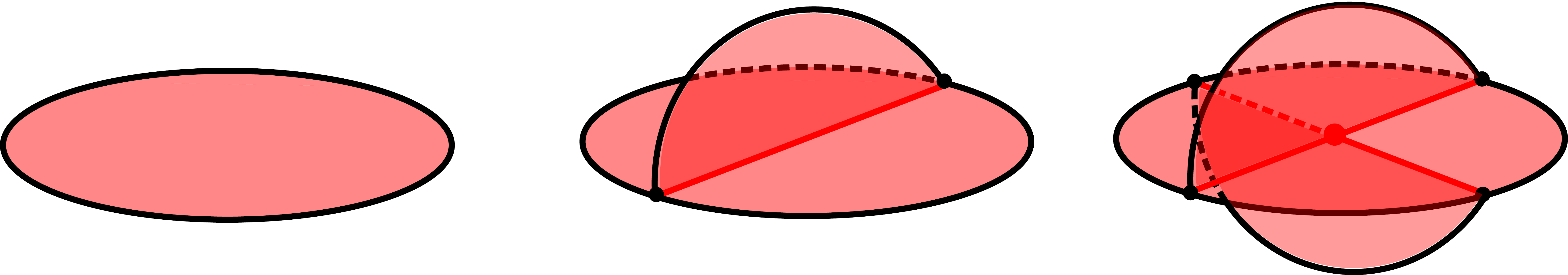}
\caption{\footnotesize The local models of the singular foam $\Sigma$: points in a 2-, 1-, and 0-stratum of $\Sigma$.}
\label{fig: Sigma_local}
\end{figure}

Second, of course the hyperbolic, half-pipe and Anti-de Sitter structures have singularities. Unfortunately here we do not only have cone singularities, and the singular locus $\Sigma\subset\mathcal M$ is slightly more complicated. In fact $\Sigma$ is a two-dimensional simplicial complex, which is a \emph{foam}, that is, it is locally homeomorphic to the cone over the skeleton of a tetrahedron.  This means that there is a zero measure subset of $\Sigma$ where more complicated singularities occur: see Figure \ref{fig: Sigma_local}. However, this is in some sense the simplest type of singularities that may occur besided cone singularities. Moreover $\Sigma$ is compact, meaning that the singularities do not enter into the cusps.

\bibliographystyle{alpha}
\bibliographystyle{ieeetr}
\bibliography{sr-bibliography}

\end{document}